\documentclass[reqno,centertags]{amsart}
\usepackage{amsmath,amsthm,amscd,amssymb}
\usepackage{latexsym}
\usepackage{hyperref}

\newcommand{\bbC}{{\mathbb{C}}}
\newcommand{\Cmplx}{{\mathbb{C}}}
\newcommand{\bbD}{{\mathbb{D}}}

\newcommand{\bbR}{{\mathbb{R}}}
\let\Reals=\bbR

\newcommand{\bbZ}{{\mathbb{Z}}}

\newcommand{\dott}{\,\cdot\,}

\newcommand{\lb}{\label}
\newcommand{\f}{\frac}

\newcommand{\ol}{\overline}
\newcommand{\ti}{\tilde  }

\newcommand{\tr}{\text{\rm{Tr}}}

\newcommand{\loc}{\text{\rm{loc}}}

\newcommand{\ran}{\text{\rm{ran}}}

\newcommand{\ess}{\text{\rm{ess}}}
\newcommand{\ac}{\text{\rm{ac}}}

\newcommand{\s}{\text{\rm{s}}}

\newcommand{\supp}{\text{\rm{supp}}}

\newcommand{\bi}{\bibitem}

\newcommand{\beq}{\begin{equation}}
\newcommand{\eeq}{\end{equation}}
\newcommand{\ba}{\begin{align}}
\newcommand{\ea}{\end{align}}
\newcommand{\veps}{\varepsilon}

\newcounter{smalllist}
\newenvironment{SL}{\begin{list}{{\rm(\roman{smalllist})}}{%
\setlength{\topsep}{0mm}\setlength{\parsep}{0mm}\setlength{\itemsep}{0mm}%
\setlength{\labelwidth}{2em}\setlength{\leftmargin}{3em}\usecounter{smalllist}%
}}{\end{list}}

\DeclareMathOperator{\Real}{Re}
\DeclareMathOperator{\Ima}{Im}
\let\Re=\undefined\DeclareMathOperator*{\Re}{Re}
\let\Im=\undefined\DeclareMathOperator*{\Im}{Im}

\DeclareMathOperator*{\slim}{s-lim}
\DeclareMathOperator*{\wlim}{w-lim}

\allowdisplaybreaks

\numberwithin{equation}{section}

\newtheorem{theorem}{Theorem}[section]

\newtheorem*{p2.1}{Proposition 2.1}
\newtheorem{proposition}[theorem]{Proposition}
\newtheorem{lemma}[theorem]{Lemma}
\newtheorem{corollary}[theorem]{Corollary}
\theoremstyle{definition}

\theoremstyle{remark}
\newtheorem*{remark}{Remark}
\newtheorem*{remarks}{Remarks}

\theoremstyle{definition}

\newcommand{\abs}[1]{\lvert#1\rvert}

\begin{document}
\title[Sum Rules and Spectral Measures]
{Sum Rules and Spectral Measures of Schr\"odinger Operators with $L^2$ Potentials}
\author[R. Killip and B. Simon]{Rowan Killip$^{1}$ and Barry Simon$^{2}$}

\thanks{$^1$ Department of Mathematics, UCLA, Los Angeles, CA 90095.
E-mail: killip@math.ucla.edu. Supported in part by NSF grant DMS-0401277 and a Sloan Foundation Fellowship.}
\thanks{$^2$ Mathematics 253-37, California Institute of Technology, Pasadena, CA 91125.
E-mail: {\nobreak bsimon@caltech.edu}. Supported in part by NSF grant DMS-0140592 and in part
by Grant No.\ 2002068 from the United States--Israel Science Foundation, Jerusalem, Israel}

\date{August 30, 2006}

\begin{abstract} Necessary and sufficient conditions are presented for a positive measure
to be the spectral measure of a half-line Schr\"odinger operator with square integrable potential.
\end{abstract}

\maketitle

\section{Introduction} \lb{s1}

In this paper, we will discuss which measures occur as the spectral measures for half-line
Schr\"odinger operators with certain decaying potentials.   Let us begin with the appropriate
definitions.

A potential $V\in L^2_\loc(\Reals^+)$ (where $\Reals^+=[0,\infty)$) is said to be limit point at infinity if
\begin{equation}
H=-\frac{d^2\ }{dx^2} + V(x)
\end{equation}
together with a Dirichlet boundary condition at the origin, $u(0)=0$, defines a selfadjoint
operator on $L^2(\Reals^+)$, without the need for a boundary condition at infinity.
This is what we will mean by a Schr\"odinger operator.  (In some
sections, we will also treat $V\in L^1_\loc$ where we feel that this generality may be of use
to others.)

The spectral theory of such operators was first described by Weyl and
subsequently refined by many others.  We will now sketch the parts of this theory that
are required to state our results; fuller treatments can be found elsewhere (e.g.,
\cite{CL,RS2,Titchmarsh}).

The name `limit point' was coined by Weyl for the following property, which is equivalent to that
given above: for all $z\in\Cmplx\setminus\Reals$ there exists a unique function $\psi\in
L^2(\Reals^+)$ so that $-\psi'' + V\psi=z\psi$ and $\psi(0)=1$.  The value of $\psi'(0)$ is denoted
$m(z)$ and is termed the (Weyl) $m$-function. It is an analytic function of $z$.
Of course, by homogeneity, one has that
\begin{equation}\label{mDefn}
m(z) = \frac{\psi'(0)}{\psi(0)}
\end{equation}
where $\psi$ is \emph{any} non-zero $L^2$ solution of $-\psi'' + V\psi=z\psi$.  This will prove the more convenient
definition.

Simple Wronskian calculations show that $m(z)$ has a positive imaginary part whenever $\Im(z)>0$.
Therefore, by the Herglotz Representation Theorem, there is a unique positive measure $d\rho$ so that
\begin{equation}\label{HRTa}
    \int \frac{d\rho(E)}{1+E^2} < \infty
\end{equation}
and
\begin{equation}\label{HRT}
m(z) = \int \ \biggl[ \frac{1}{E-z} - \frac{E}{1+E^2} \biggr] \, d\rho(E) + \Re m(i).
\end{equation}
Uniqueness follows from the fact that
\begin{equation} \label{RBV}
d\rho(E) =\wlim_{\varepsilon\downarrow 0}\, \tfrac{1}{\pi}\, \Im m(E+i\varepsilon)\, dE.
\end{equation}
Moreover, the boundary values $m(E+i0)$ exist almost everywhere and are equal to the Radon-Nikodym
derivative, $\frac{d\mu}{dE}$.

At first sight, \eqref{HRT} does not permit us to recover $\Im m$ from $d\rho$ without first knowing
$\Re m(i)$.  Actually, it can be recovered from the asymptotic \cite{Atk,GSAnn}
\begin{equation} \lb{1.10}
m(z) =\sqrt{-z} + o(1)
\end{equation}
that holds as $|z|\to\infty$ along rays at a small angle to the negative real axis.  (If the support
of $d\rho$ is bounded from below, it holds as $z\to-\infty$.)

Just as $V$ determines $d\rho$, so $\rho$ (or $m(z)$) determines $V$.  This is a famous result of
Gel'fand--Levitan \cite{GL55,LevBk,LevSar}; see also Remling \cite{Rem02} and Simon \cite{Sim271}.

As we have described, each potential gives rise to a spectral measure, which also
determines $V$.  Our main goal in this paper is to give necessary and sufficient conditions in
terms of $\rho$ for $V\in L^2(\Reals^+)$.  Our model here is a result we proved  recently
\cite{KS} for Jacobi matrices, the discrete analog of Schr\"odinger operators.  (Other pre-cursors
will be discussed later).  To properly frame our result, we recall
briefly the Jacobi matrix result. A Jacobi matrix is a semi-infinite tridiagonal matrix
\begin{equation} \lb{1.11}
J= \begin{pmatrix}
b_1 & a_1 & 0   & \cdots \\
a_1 & b_2 & a_2 & \vphantom{\ddots} \\
0   & a_2 & b_3 & \ddots \\
\vdots &   & \ddots & \ddots
\end{pmatrix}
\end{equation}
viewed as an operator on $\ell^2 (\bbZ_+)=\ell^2 (\{1,2,3,\dots\})$. The spectral measure
here is defined by
\begin{equation} \lb{1.12}
m(E) =\langle \delta_1, (J-E)^{-1}\delta_1\rangle =\int\f{d\mu(x)}{x-E}.
\end{equation}
We write $J_0$ for the Jacobi matrix with $b_n\equiv 0$, $a_n\equiv 1$.

Our earlier result is:

\begin{theorem}[\cite{KS}]\lb{T1.1} $J-J_0$ is Hilbert-Schmidt, that is,
\begin{equation} \lb{1.13}
\sum_{n=1}^\infty (a_n -1)^2 + b_n^2 <\infty,
\end{equation}
if and only if the spectral measure $d\mu$ obeys
\begin{SL}
\item {\rm{(Blumenthal--Weyl)}} $\supp (d\mu) = [-2,2] \cup
\{E_j^+\}_{j=1}^{N_+} \cup \{E_j^-\}_{j=1}^{N_-}$ with $E_1^+ > E_2^+ >
\cdots > 2$ and $E_1^- < E_2^- < \cdots < -2$ with $\lim_{j\to\infty} E_j^\pm
=\pm2$ if $N_\pm =\infty$.

\item {\rm{(Normalization)}} $\mu$ is a probability measure.

\item {\rm{(Lieb--Thirring Bound)}}
\begin{equation}\lb{1.12x}
\sum_{\pm, j} (\abs{E_j^\pm} -2)^{3/2}  <\infty
\end{equation}
\item {\rm{(Quasi-Szeg\H{o} Condition)}}  Let $d\mu_{\ac}(E)=
f(E)\, dE$. Then
\begin{equation} \lb{1.14}
\int_{-2}^2 \log (f(E)) \sqrt{4-E^2}\, dE >-\infty
\end{equation}
\end{SL}
\end{theorem}

We have changed the ordering of the conditions relative to \cite{KS} in order to
facilitate comparison with Theorem~\ref{T1.2} below.

To state the result for Schr\"odinger operators, we need some further preliminaries.
Let $d\rho_0$ be the free spectral measure (i.e., for $V=0$), it is
\begin{equation} \lb{1.19}
d\rho_0 (E) =\pi^{-1} \chi_{[0,\infty)}(E) \sqrt{E}\, dE.
\end{equation}
For any positive measure $\rho$, we define a signed measure $d\nu$ on $(1,\infty)$ by
\begin{equation} \lb{1.17a}
\f{2}{\pi} \int f(k^2) k\,d\nu(k) = \int f(E) [d\rho(E)-d\rho_0 (E)].
\end{equation}
Notice that $d\nu$ is parameterized by momentum, $k$, rather than energy,
$E=k^2$.  This is actually the natural independent variable for what follows (in \cite{KS}
we used $z$ defined by $E=z+z^{-1}$).   We will write $w$ for the $m$-function in
terms of $k$:
\begin{equation}
w(k) = m(k^2).
\end{equation}
With this notation,
\begin{equation} \lb{1.17b}
\f{d\nu}{dk} = \Im[ w(k+i0) ] - k
\end{equation}
at a.e.\ point $k\in(1,\infty)$.  Here $\frac{d\nu}{dk}$ is the Radon-Nikodym
derivative of the a.c.\ part of $\nu$, which hay have a singular part as well.

We will eventually prove (see Section~\ref{s9}) that if $V\in L^2$, then
\begin{equation} \lb{1.17c}
\bigl\| \nu \bigr\|_{\ell^2(M)}^2 := \sum \bigl[ \abs{\nu} ([n,n+1]) \bigr]^2 <\infty
\end{equation}
and as a partial converse, if $d\rho$ is supported in $[-a,\infty)$ and \eqref{1.17c} holds,
then $d\rho$ is the spectral measure for a potential $V\in L_\loc^2$.

We also need to introduce the long- and short-range parts of the Hardy--Littlewood
maximal function \cite{HL,Rudin},
\begin{align}
(M\nu)(x) &= \sup_{L>0}\, \f{1}{2L}\, \abs{\nu} ([x-L,x+L]) \lb{1.17d} \\
(M_s\nu)(x) &= \sup_{0<L\leq 1} \, \f{1}{2L}\, \abs{\nu} ([x-L,x+L]) \lb{1.17e} \\
(M_l\nu)(x) &= \sup_{1\leq L}\, \f{1}{2L}\, \abs{\nu} ([x-L,x+L]). \lb{1.17f}
\end{align}

The main theorem of this paper is

\begin{theorem}\lb{T1.2} A positive measure $d\rho$ on $\Reals$ is the spectral measure associated
to a $V\in L^2(\Reals^+)$ if and only if
\begin{SL}
\item {\rm{(Weyl)}} $\supp (d\rho) = [0,\infty) \cup \{E_j\}_{j=1}^N$
with $E_1 < E_2 < \cdots <0$ and  $E_j \to 0$ if $N=\infty$.

\item {\rm{(Normalization)}}
\begin{equation} \lb{1.17g}
\int \log \biggl[ 1+ \biggl( \f{M_s\nu (k)}{k}\biggr)^2\biggr] k^2\, dk <\infty
\end{equation}

\item {\rm{(Lieb--Thirring)}}
\begin{equation}\lb{1.20x}
\sum_j \, \abs{E_j}^{3/2}  <\infty
\end{equation}

\item {\rm{(Quasi-Szeg\H{o})}}
\begin{equation} \lb{1.17h}
\int_0^\infty \log \biggl[ \f{1}{4}\, \f{d\rho}{d\rho_0} + \f12 + \f14\, \f{d\rho_0}{d\rho} \biggr]
    \sqrt{E}\, dE <\infty
\end{equation}
\end{SL}
\end{theorem}

\begin{remarks} 1. It may be surprising that we have replaced the innocuous normalization
condition in the Jacobi case by \eqref{1.17g}. The reason is the following: $\mu(\bbR)=1$
in the Jacobi case is the condition that $\mu$ is the spectral measure of some
Jacobi matrix. In this theorem, we do not presume a priori that $d\rho$ is spectral
measure. We will eventually see that \eqref{1.17g} implies that $\rho$ is the spectral
measure of an $L_\loc^2$ potential. Indeed, \eqref{1.17g} has additional information
that we will need to control high energy pieces.

2. The name of condition (i) was chosen because the fact that it is implied by $V\in L^2$
is an immediate consequence of Weyl's Theorem on the invariance of the essential spectrum under
(relatively) compact perturbations.

3. Bounds on sums of powers of eigenvalues in terms of the $L^p$ norm of the potential are
usually referred to as Lieb--Thirring Inequalities in deference to their exhaustive work on
this question, \cite{LT}. However, the particular case that appears in Theorem~\ref{T1.2}
was first observed by Gardner, Greene, Kruskal, and Miura; see \cite[p.~115]{GGKM}.

4. The argument of $\log$ in \eqref{1.17h} has the form
\begin{equation}\label{1.?}
\tfrac14\, \lambda + \tfrac12 + \tfrac14\,\lambda^{-1} =
[\tfrac12\, (\lambda+\lambda^{-1})]^2 \geq 1
\end{equation}
so the integrand is nonnegative. This is significantly different from \eqref{1.14} where
the integrand can have both signs, and the finiteness of the measure implies one sign
is automatically finite so we do not have to worry about oscillations. In our case,
an oscillating integrand would present severe difficulties because spectral measures are not
finite.

5. Theorem~\ref{T1.2} implies that if $V\in L^2$, then $\sigma_\ac (H)=[0,\infty)$.
This is a result of Deift--Killip \cite{DeiftK}.
\end{remarks}

We will prove Theorem~\ref{T1.2} in two parts. First, we prove an equivalence of
$V\in L^2$ and a set of conditions that has an unsatisfactory element. Then we will
show the conditions of Theorem~\ref{T1.2} are equivalent to those of this
Theorem~\ref{T1.3}.  Both portions are lengthy.

The intermediate theorem requires one further object.  Let
\begin{equation} \lb{1.18}
    F(q) = \pi^{-1/2} \int_{p\geq 1} p^{-1} e^{-(q-p)^2} d\nu (p).
\end{equation}
By \eqref{HRTa}, this integral is absolutely convergent.

\begin{theorem}\lb{T1.3} A positive measure $d\rho$ on $\Reals$ is the spectral measure associated
to a $V\in L^2(\Reals^+)$ if and only if
\begin{SL}
\item {\rm{(Weyl)}} $\supp (d\rho) = [0,\infty) \cup \{E_j\}_{j=1}^N$ with
$E_1 < E_2 < \cdots <0$ and $E_j\to0$ if $N=\infty$.

\item {\rm{(Local Solubility)}}
\begin{equation} \lb{1.19x}
\int_0^\infty \, \abs{F(q)}^2 \, dq <\infty
\end{equation}

\item {\rm{(Lieb--Thirring)}}
\begin{equation}\lb{1.20}
\sum_j \, \abs{E_j}^{3/2}  <\infty
\end{equation}

\item {\rm{(Strong Quasi-Szeg\H{o})}}
\begin{equation} \lb{1.21}
\int \log \biggl[ \f{\abs{w(k+i0) + ik}^2}{4k\Im w(k+i0)}\biggr] k^2\, dk <\infty
\end{equation}
\end{SL}
\end{theorem}

\begin{remarks} 1. Notice that
\begin{align*}
\abs{w(k+i0)+ik}^2 &\geq \abs{\Im w(k+i0) + k}^2 \\
&\geq \abs{\Im w(k+i0) + k}^2 - \abs{\Im w(k+i0) -k}^2 \\
&= 4k \Im w(k+i0)
\end{align*}
so the argument in the $\log$ in \eqref{1.21} is at least $1$ and the integrand is
strictly positive, so the integral either converges or is $+\infty$.

2. There is a significant difference between \eqref{1.21} and \eqref{1.14}. Since
\eqref{1.21} involves $M$ and not just $\Im w$, the singular part of $d\rho$ enters
in both (ii) and (iv). Still, as we shall see, the restriction on $d\rho_\s$ is mild.

3. The occurrence of $w$ in \eqref{1.21} means that if one starts with $d\rho$,
it is difficult to check this condition---one first has to calculate the Hilbert transform
(conjugate function) of $d\rho$.  Considering the example
$$
d\rho = d\rho_0 + \sum_{j=1}^\infty c_j \delta(E-j^2) dE
$$
shows that \eqref{1.19x} is not strong enough to allow the replacement of \eqref{1.21} by
the weaker Quasi-Szeg\H{o} condition \eqref{1.17h}.  Specifically, \eqref{1.19x} only requires $c_j\in\ell^2$
while \eqref{1.21} implies $c_j\in\ell^1$.  The relation of the two Quasi-Szeg\H{o} conditions
and their connection to $\Re w(k)$ is discussed in Section~\ref{s8}.

The advantage of the maximal function is that it involves no cancellation; we see plainly that
\eqref{1.17g} is a statement about the size of $d\rho-d\rho_0$.

4. The name `Local Solubility' comes from the fact that this condition (plus the fact that
support of $d\rho$ is bounded from below) guarantees that $d\rho$ is the spectral measure
for some $L^2_\loc$ potential.  See Section~\ref{s6}.

5. We will prove Theorem~\ref{T1.3} in Section~\ref{s7}, and then use it to prove
Theorem~\ref{T1.2} in Sections~\ref{s8}--\ref{s11}.
\end{remarks}

There are significant differences from Theorem~\ref{T1.1}, both in the form of
Theorem~\ref{T1.3} and its proof. Understanding the difficulties that led to these
differences is illuminating. To understand the issues, we recall that Theorem~\ref{T1.1} was
proven by showing a general sum rule, dubbed the $P_2$ sum rule:
\begin{equation} \lb{1.22}
Q(d\mu) + \sum_{\pm,j} F(E_j^\pm) =\tfrac14\, \sum_j b_j^2 + \tfrac12\, \sum_j G(a_j)
\end{equation}
where
\begin{align}
Q(d\mu) &=\f{1}{4\pi} \int_{-2}^2 \log \biggl( \f{\sqrt{4-E^2}}{2\Ima m(E+i0)}\biggr)
\sqrt{4-E^2}\, dE \lb{1.23} \\
G(a) &= a^2 -1 -\log \abs{a}^2 \lb{1.24} \\
F(\beta + \beta^{-1}) &= \tfrac14\, [\beta^2 -\beta^{-2} -\log \beta^4 ], \quad |\beta|>1. \lb{1.25}
\end{align}
To prove Theorem~\ref{T1.1}, one
proves \eqref{1.22} is always true (both sides may be infinite) and then notes
$Q(d\mu)<\infty$ if and only if \eqref{1.14} holds, $F(E_j^\pm)=(\abs{E_j^\pm}-2)^{3/2}
+ (O(\abs{E_j^\pm}-2)^2)$, and $G(a)=2(a-1)^2 + O((a-1)^3)$, so $\sum_{\pm,j}
F(E_j^\pm)<\infty$ if and only if \eqref{1.12x} holds and the right side of
\eqref{1.13} is finite if and only if \eqref{1.11} holds.

For nice $J$'s, (e.g., $b_n=0$ and $a_n=1$ for $n$ large), \eqref{1.22} is a
combination of two sum rules of Case \cite{Case1,Case2}. For general $J$'s,
it is proven in Killip--Simon \cite{KS} with later simplifications of parts of
the proof in \cite{Sim288,SZ}.

The difficulties in extending this strategy in the continuous case were several:
\begin{SL}
\item The translation of the normalization condition $\mu(\bbR)=1$ is
not clear.  We needed a condition that guaranteed $d\rho$ is the spectral measure associated
to a reasonable $V$\!, preferably belonging to $L^2_\loc$.  We sought to express this in terms
of the divergence of $\rho(-\infty,R)$ as $R\to\infty$. As it turned out, the $A$-function approach
to the inverse spectral problem, \cite{GSAnn,Sim271}, leads quickly and conveniently to the condition
\eqref{1.19x}, which is perfect for us.

\item The natural half-line sum rules in the Schr\"odinger case
invariably lead to terms involving $V(0)$ or worse still, $V'(0)$.  This is clearly unacceptable
for one seeking $V\in L^2$ conditions.

\item The half-line sum rules also lead to terms that, like \eqref{1.14},
have an integrand that has a variable sign. In \eqref{1.14}, the fact that
$\int f(E)\, dE \leq 1$ implies uniform control on $\int \log_+ (f(E)) \sqrt{4-E^2}\,
dE$ and so the terms of the `wrong' sign (where $f(E)>1$) present no problem.
But in the whole-line case where $\rho (\bbR)=\infty$, terms of opposite signs
could involve difficult to control cancellations.
\end{SL}

The resolution of difficulties (ii) and (iii) was to fall back to the whole-line sum rule used in
\cite{DeiftK}.  The penalty is that the Strong Quasi-Szeg\H{o} condition, \eqref{1.21}, little resembles
the Quasi-Szeg\H{o} condition of our earlier theorem, \eqref{1.14}.  It is this disappointment that led
us to push on and find Theorem~\ref{T1.2}.

Whole-line sum rules date to the original inverse-scattering solution of the KdV equation, \cite{GGKM}.
Consider the operator
\begin{equation} \lb{1.26}
L_0 =-\f{d^2}{dx^2} +\chi_{(0,\infty)} (x) V(x)
\end{equation}
acting on $L^2(\Reals)$ with eigenvalues $E_j^{(0)}$. The well-known sum rule is, \cite{DeiftK,GGKM,ZF},
\begin{equation} \lb{1.27}
\tfrac18\, \int_0^\infty V(x)^2\, dx = \tfrac23\, \sum_j (E_j^{(0)})^{3/2} + Q
\end{equation}
where
\begin{equation} \lb{1.28}
Q =\f{1}{\pi} \int_0^\infty \log
\biggl[ \f{\abs{w(k+i0)+ik}^2}{4k\Im w(k+i0)}\biggr] k^2\, dk.
\end{equation}

As in \cite{KS}, we will need to prove it in much greater generality than
was known previously. Essentially, assuming that $w$ is the $m$-function of an
$L_{\loc}^2$ potential $V$\!, we will prove \eqref{1.27} always holds although
both sides may be infinite.

If one notes that the half-line and whole-line eigenvalues interlace,
\begin{equation} \lb{1.29}
E_j^{(0)}\leq E_j \leq E_{j+1}^{(0)},
\end{equation}
it is clear that \eqref{1.27} proves Theorem~\ref{T1.2}.

As was the case \cite{KS,Sim288,SZ}, the key to the proof of \eqref{1.27} is a
`step-by-step' sum rule, that is, a result that, in essence, is the difference
of \eqref{1.27} for $L_0$ and for
\begin{equation} \lb{1.30}
L_t =-\f{d^2}{dx^2} +\chi_{(t,\infty)}(x) V(x)
\end{equation}
and which always holds.  A second important ingredient is the semicontinuity of $Q$.

In Section~\ref{s2}, we will discuss a relative Wronskian which is the analogue
of the product of $m$-functions used implicitly in \cite{KS,SZ} and explicitly in
\cite{Sim288} to prove a multi-step sum rule. In Section~\ref{s3}, as
an aside, we will re-express this relative Wronskian as a perturbation determinant.
In Section~\ref{s4}, we prove the step-by-step sum rule. In Section~\ref{s5}, we
prove lower semicontinuity of the quasi-Szeg\H{o} term. In Section~\ref{s6}, we discuss
\eqref{1.19x} and, in particular, show it implies $\mu$ is the spectral measure of
a locally $L^2$ potential. Section~\ref{s7} completes the proof of \eqref{1.27} and
of Theorem~\ref{T1.3}. Sections~\ref{s8}--\ref{s11} prove Theorem~\ref{T1.2} given
Theorem~\ref{T1.3}.

The earliest theorem of the type presented here is Verblunsky's form \cite{V36} of
Szeg\H{o}'s theorem \cite{OPUC1,Szb}.  Let us elaborate.  The orthogonal polynomials associated to a
measure on the unit circle obey a recurrence and the coefficients that appear in this
recurrence are known as the Verblunsky coefficients.  The result just mentioned says
that the Verblunsky coefficients are square summable if and only if the logarithm of
the density of the a.c.\ part of the measure is integrable.  In fact, there is a sum
rule relating these quantities.

One of the more interesting spectral consequences of Szego's theorem is the construction by
Totik \cite{Totik} (see also Simon \cite{OPUC1}) that given any measure supported on the
circle, there is an equivalent measure whose recursion coefficients lie in
all $\ell^p$ ($p>2$).  We expect that the results and techniques of the current
paper will provide tools allowing one to carry this result over to
Schr\"odinger operators (although it seems likely that $\ell^p$ will be replaced by
$\ell^p(L^2)$ rather than by $L^p$.

Kre{\u\i}n  systems give a continuum analogue for orthogonal polynomials on the unit
circle.  The corresponding version of Szeg\H{o}'s Theorem can be found in \cite{Krein}; though
for proofs, see \cite{Stakh}.  Using a continuum analogue of the Geronimus relations,
Kre{\u\i}n's Theorem gives results for potentials of the form $V(x)=a(x)^2 \pm a'(x)$
with $a\in L^2$.  Note that the operators associated to such potentials are automatically
positive---there are no bound states.  For a further discussion of the application of
Kre{\u\i}n  systems to Schr\"odginer operators, see \cite{DenJMAA,DenCMP,DenJDE}.

More recently, Sylvester and Winebrenner \cite{SW} studied the scattering for the
Helmholtz equation on a half-line and obtained necessary and sufficient conditions (in
terms the reflection coefficient) for square integrability of the derivative of the wave speed.
Applying appropriate Liouville transformations connects this work to the study of Schr\"odinger
operators with potentials $V(x)=a(x)^2 \pm a'(x)$, just as for Kre{\u\i}n  systems.  Our methods
parallel their work in places, particularly with regard to the semicontinuity properties of
$Q$ discussed in Section~\ref{s5}.  However, dealing with bound states adds to the complexity of
our case.

As mentioned earlier, \cite{DeiftK} proved that
$\sigma_{\ac}(H)=[0,\infty)$ for $H=-\f{d^2}{dx^2}+V$ with $V\in L^2$.  Earlier
work by Christ, Kiselev, and Remling, \cite{Kis96,CK,Rem98}, settled the case
$V(x) \leq C (1+\abs{x}^2)^{-\alpha}$ for $\alpha>\f12$ by entirely different means.
The most recent development in this direction is the use of sum rules by Rybkin, \cite{Rybppt},
to prove $\sigma_{\ac}(H)=[0,\infty)$ for potentials of the form $V=f+g'$ with $f,g\in L^2$.

\medskip
\noindent
\textit{Acknowledgements:}
We wish to thank Wilhelm Schlag, Terence Tao, and Christoph Thiele for various pointers
on the harmonic analysis literature.  We would also like thank Christian Remling for some
insightful comments.

\section{The Relative Wronskian} \lb{s2}

In this section, we will consider $V\in L_\loc^1 (\Reals^+)$ for which the operator
\begin{equation} \lb{2.1}
H=-\f{d^2}{dx^2}+V
\end{equation}
with boundary condition $u(0)=0$ is essentially selfadjoint and has
$\sigma_\ess (H)\subset [0,\infty)$. As noted in the Introduction, for any $k\in\bbC_+$
with $k^2\not\in\sigma(H)$, there is unique solution $\psi_+ (x,k)$ of
\begin{equation}\label{R2.1}
-\psi'' + V\psi = k^2 \psi
\end{equation}
which is $L^2$ at $+\infty$ and $\psi_+(0)=1$.  By the above assumption on $\sigma_\ess$,
this extends to a meromorphic function of $k$ in $\bbC_+$ with poles exactly at the negative
eigenvalues of $H$.  Moreover, the poles are simple.

Let us define
\begin{equation} \lb{2.2}
W(x,k) =e^{-ikx}\psi'_+ (x,k) + ike^{-ikx}\psi_+ (x,k).
\end{equation}
$W$ is the Wronskian of $\psi_+ (x,k)$ and $\psi_-^{(0)}(x,k) \equiv e^{-ikx}$, which is the solution of
\begin{equation} \lb{2.2a}
-\psi'' =k^2 \psi
\end{equation}
that is $L^2$ at $-\infty$ (recall $\Im k>0$).   Note that $W(x,k)$ is a meromorphic function of $k$,
an absolutely continuous function of $x$, and is easily seen to obey
\begin{equation} \lb{2.3}
\f{\partial}{\partial x}\, W(x,k) =e^{-ikx} \psi_+ (x,k) V(x).
\end{equation}

The zeros of $k\mapsto W(x_0,k)$ are precisely those points where one can find a
$c\in\Cmplx$ for which
\begin{equation} \lb{2.4}
u(x) = \begin{cases} \psi_+ (x,k), & x\geq x_0 \\
ce^{-ikx}, & x\leq x_0 \end{cases}
\end{equation}
is a $C^1$ function, that is, $W(x_0,k)=0$ if and only if $k^2$ is an eigenvalue
of the operator $L_t$ of \eqref{1.27} with $t=x_0$. In particular, all
zeros lie on the imaginary axis: $k=i\kappa$ with $\kappa >0$.

We will use $\kappa_1(x) >\kappa_2(x) >\cdots$ to indicate the zeros
of $W(x,k)$ so that $-\kappa_j(x)^2$ are the negative eigenvalues of $L_x$.

We define the \textit{relative Wronskian} by
\begin{equation} \lb{2.5}
a_x (k) =\f{W(x,k)}{W(0,k)}.
\end{equation}
For each $x$, it is a meromorphic function of $k$. Like the $m$-function---and unlike $W(x,\cdot)$---it is
independent of the normalization $\psi_+(0,k)=1$.  By the above, we have

\begin{proposition} \lb{P2.1} The poles of $a_x(k)$ are simple and lie at those points $k=-i\kappa_j(0)$ for which
$-\kappa_j(0)^2$ is an eigenvalue of $L_0$.   The zeros are also simple and lie on the set $k=-i\kappa_j(x)$ where
$-\kappa_j(x)^2$ are the eigenvalues of $L_x$.  {\rm{(}}In the event that a point lies in both sets, there is neither
a pole nor a zero---they cancel one another.{\rm{)}}
\end{proposition}

Next, we note that

\begin{proposition}\lb{P2.2} For $\Re k\neq 0$, $a_x (k)$ is absolutely
continuous in $x$. Moreover, $\log[a_x (k)]$, defined so that it is continuous
in $x$ with $\log (a_{x=0}(k))=0$, obeys
\begin{equation} \lb{2.6}
\f{d}{dx}\, \log[a_x(k)] = V(x) (ik+w(k;x))^{-1}
\end{equation}
where
\begin{equation}\label{E:w}
w(k,x) = \frac{\psi_+'(x,k)}{\psi_+(x,k)}
\end{equation}
is the $m$-function associated to the operator $L_x$ restricted to $[x,\infty)$.
\end{proposition}

\begin{proof} As a ratio of nonvanishing absolutely continuous functions, $a_x(k)$
is absolutely continuous, and then so is its log. By \eqref{2.3},
\[
\f{d}{dx}\, \log[a_x(k)] = \f{e^{-ikx} \psi_+ (x,k) V(x)}{W(x,k)}.
\]
As
\[
\f{W(x,k)}{e^{-ixk} \psi_+(x,k)} = \f{\psi'_+ (x,k)}{\psi_+ (x,k)} + ik = w(k;x)+ik,
\]
\eqref{2.6} is immediate.
\end{proof}

\begin{proposition} \lb{P2.3}
\begin{SL}
\item[{\rm{(a)}}] For $k\in\bbC_+$ with $\Real k\neq 0$,
\begin{equation} \lb{2.6a}
\bigl|  \log[a_x(k)] \bigr| \leq \abs{\Real k}^{-1} \int_0^x \abs{V(y)}\, dy.
\end{equation}
\item[{\rm{(b)}}] Fix $K>0$. Then there exist $R>0$ and $C$ so, for all $k$ in
$\bbC_+$ with $\abs{k} >R$ and all $x$ in $[0,K]$,
\begin{equation} \lb{2.6b}
\bigl|  \log[a_x(k)] \bigr| \leq C\abs{k}^{-1}.
\end{equation}
\end{SL}
\end{proposition}

\begin{proof} (a) If $\Real k>0$ and $\Ima k>0$, then $\Ima w>0$ and so $\abs{ik+w(k;x)}^{-1}
\leq \abs{\Real k}^{-1}$.  Thus \eqref{2.6a} follows from \eqref{2.6}.

\smallskip
(b) By \cite{GSAnn}, uniformly for $x\in [0,K]$, $w(k;x)-ik\to 0$ as $\abs{k}\to\infty$
with $\arg (-ik) \leq \f{\pi}{4}$. This plus \eqref{2.6} implies that \eqref{2.6b} holds
uniformly in $x\in (0,K)$ and $\abs{\arg (-ik)}\leq \f{\pi}{4}$. By
\eqref{2.6a}, it holds for $\arg k\in (0, \f{\pi}{4})\cup (\f{3\pi}{4},\pi)$.
\end{proof}

\begin{proposition} \lb{P2.4} Suppose that for some $k_0\in (0,\infty)$,
$\lim_{\veps\downarrow 0} w(k_0 + i\veps; 0)\equiv w(k_0 +i0; 0)$ exists and
$\Ima w(k_0 +i0;0)\in (0,\infty)$. Then for all $x$, $\lim_{\veps\downarrow 0}
w(k_0 +i\veps;x) \equiv w(k_0 +i0; x)$ exists and $\lim_{\veps\downarrow 0} a_x
(k_0 +i\veps) \equiv a_x (k_0 +i0)$ exists. Moreover,
\begin{equation} \lb{2.7}
\abs{a_x (k_0 + i0)}^2 = \f{T(k_0,0)}{T(k_0,x)}
\end{equation}
where
\begin{equation} \lb{2.8}
T(k_0,x) = \f{4k_0 \Ima w(k_0 +i0;x)}{\abs{w(k_0 +i0;x)+ik_0}^2}.
\end{equation}
\end{proposition}

\begin{proof} As $\Ima w(k_0 +i0;0)\in (0,\infty)$, we may also take the vertical limit
$\psi_+(x,k_0+i0)$; indeed, this is just the solution to \eqref{R2.1} with $\psi(0)=1$,
$\psi'(0)=\Ima w(k_0 +i0;0)$, and $k=k_0$.

As $\Ima w(k+i0)>0$, $\psi_+ (x,k_0 +i0)$ is not a complex multiple of a real-valued solution and so
cannot have any zeros. Thus $\lim_{\veps\downarrow 0} w(k_0 + i\veps;x)$ exists. Similarly,
by \eqref{2.2}, $W(x,k_0 +i\veps)$ has a limit and
\begin{equation} \lb{2.9}
\abs{W(x, k_0 +i0)}=\abs{\psi_+ (x, k_0 +i0)}\, \abs{w(k_0 +i0;x)+ik}.
\end{equation}

As $\psi_+ (\dott, k_0 +i0)$ and $\ol{\psi_+ (\dott, k_0 +i0)}$ obey the
same equation, their Wronskian is a constant (in $x$), that is
\begin{equation} \lb{2.10}
\abs{\psi_+ (x, k_0 +i0)}^2 \Ima w(x;k_0 +i0)=C_{k_0}.
\end{equation}
The definition \eqref{2.8} and \eqref{2.9}, \eqref{2.10} imply
\begin{equation} \lb{2.11}
\abs{W(x,k_0 +i0)}^2 = \f{4C_{k_0} k_0}{T(k_0,x)}.
\end{equation}
Thus \eqref{2.5} implies \eqref{2.7}.
\end{proof}

We write the letter $T$ in \eqref{2.8} because, as we will see, it represents the transmission
probability of stationary scattering theory.

Our final result in this section is

\begin{proposition} \lb{P2.5} Let $V\in L_\loc^2 ([0,\infty))$ and suppose
$\sigma_\ess (H)\subset [0,\infty)$. Then as $\kappa\to\infty$ {\rm{(}}real
$\kappa${\rm{)}},
\begin{equation} \lb{2.12}
\log [a_x (i\kappa)] = -\f{1}{2\kappa} \int_0^x V(y)\, dy + \f{1}{8\kappa^3}
\int_0^x V(y)^2\, dy+ o(\kappa^{-3})
\end{equation}
with an error uniform in $x$ for $x\in [0,K]$ for any $K$.
\end{proposition}

\begin{proof} By \cite{GSAnn,Sim271},
\begin{equation}\lb{2.13}
w(i\kappa,x) = -\kappa -\int_0^1 V(x+y) e^{-2\kappa y}\, dy + o(\kappa^{-1})
\end{equation}
with the $o(\kappa^{-1})$ uniform in $x$ for $x\in [0,K]$. Notice that the
integral in \eqref{2.13} is $o(\kappa^{-1/2})$ since $V\in L_\loc^2$. Thus
\begin{equation} \lb{2.14}
[w(i\kappa,x)-\kappa]^{-1} = (-2\kappa)^{-1} + (2\kappa)^{-2} \int_0^1 V(x+y)
e^{-2\kappa y}\, dy + o(\kappa^{-3}).
\end{equation}
To get this, note that one error term is $O(\kappa^{-3}) o(\kappa^{-1})$ and
by the fact that the integral in \eqref{2.13} is a priori
$o(\kappa^{-1/2})$, the other is $O(\kappa^{-3}) O(\kappa^{-1/2})^2$.

Thus, by integrating \eqref{2.6}, the proposition will follow once we show
\begin{equation} \lb{2.15}
\lim_{\kappa\to\infty} \, \kappa \int_0^x V(y) \int_0^1 V(y+s) e^{-2\kappa s} \, ds
= \tfrac12 \int_0^x V(y)^2 \, dy
\end{equation}
for all $V\in L^2 (0, x+1)$. To prove this, note first that it is
trivial if $V$ is continuous since then, $\kappa \int_0^1 V(y+s) e^{-2\kappa s}\,
ds =\f12 [V(y) +o(1)]$ for each $y$ uniformly in $y$ in $[0,x]$. Moreover,
\begin{equation} \lb{2.16}
\begin{split}
\biggl| \, \int_0^x f(y)\biggl( \int_0^1 & g(y+s) \kappa e^{-2\kappa s}\,
ds\biggr) dy \biggr| \\
& \leq \tfrac12 \biggl( \int_0^x \abs{f(y)}^2 \, dy \biggr)^{1/2}
\biggl( \int_0^{x+1} \abs{g(y)}^2 \, dy \biggr)^{1/2}
\end{split}
\end{equation}
so an approximation theorem goes from $V$ continuous to general $V$ in
$L^2 (0,x+1)$.

To prove \eqref{2.16}, use the Schwartz inequality and
\[
\biggl\| \int_0^1 g(\dott +s) \kappa e^{-2\kappa s}\, ds \biggr\|_2 \leq \int_0^1
\kappa e^{-2\kappa s} \|g(\dott +s)\|_2 \, ds \leq \tfrac12\,
\biggl( \int_0^{x+1} \abs{g(y)}^2\, dy \biggr)^{1/2}
\]
where $\|\cdot\|_2$ is $L^2 (0,x)$ norm.
\end{proof}

\section[Perturbation Determinants]{Perturbation Determinants: An Aside} \lb{s3}

In this section, we provide an alternate definition of $a_x(k)$ which we could have
used (and, indeed, initially did use) to define and prove the basic properties of
this function. The definition as a perturbation determinant makes the similarity
to the Jacobi matrix theory stronger. Expressions of suitable Wronskians as Fredholm
determinants go back to Jost and Pais \cite{JP}.   We will not use this alternate
definition again in this paper, but felt it is suggestive and should be useful for
other purposes.

We will write $\mathfrak{I}_1$ for the space of trace-class operators with the usual
norm: $\|A\|=\tr(|A|)$.

We need one preliminary:

\begin{proposition}\lb{P3.1} Let $V$ be in $L_\loc^1 ((0,\infty))$ and consider $L = -\f{d^2}{dx^2}
+V$ on $L^2(\Reals)$ {\rm{(}}with a boundary condition at infinity if
$V$ is limit circle there{\rm{)}}.  Fix $0<K<\infty$ and view $L^2([0,K])$ as functions
{\rm{(}}and multiplication operators{\rm{)}} on all of \/$\Reals$ that happen to
vanish outside this interval.

Given $z\in\bbC_+$, the mapping
$f\mapsto f(L-z)^{-1}$  is continuous and differentiable from $L^2(0,K)$ into the trace class operators.
\end{proposition}

\begin{proof} Let $L_D$ be the operator with a Dirichlet boundary condition added
at $x=K$, that is, $L_D =L_D^- \oplus L_D^+$ with $L_D^-$ on $L^2 (-\infty, K)$ with
$u(K)=0$ boundary conditions and $L_D^+$ on $L^2 (K,\infty)$ with $u(K)=0$ boundary
conditions and the same boundary condition at infinity as $L$.

Let $u_\pm$ solve $-u'' +Vu =zu$ with $u_-$ square-integrable at $-\infty$ and $u_+$, $L^2$ at
$+\infty$ (or, obeying $H$'s boundary condition at infinity if $V$ is limit circle).
Let $\varphi$ be given by
\begin{equation} \lb{3.1}
\varphi(x) =\begin{cases} u_+ (K) u_-(x), & x\leq K \\
u_- (K) u_+(x), & x\geq K \end{cases}
\end{equation}
and normalize $u_-$ so that $W(u_+, u_-)=1$. Then, standard formulae for Green's functions
\cite{CL} show that with $G(x,y)$, the integral kernel of $(L-z)^{-1}$ and $G_D(x,y)$ that of $(L_D-z)^{-1}$,
\begin{equation} \lb{3.2}
G(x,y) -G_D(x,y) = (u_+(K) u_-(K))^{-1} \varphi (x) \varphi(y).
\end{equation}

Since $\varphi$ is bounded on $[0,K]$, $f\varphi\in L^2$ and so
\[
f[(L-z)^{-1} - (L_D-z)^{-1}]
\]
is a bounded rank one operator, and so trace class. Thus it suffices to prove
$f(L_D-z)^{-1} = f(L_D^- -z)^{-1} \oplus 0$ is trace class, and so that
$f(L_D^- -z)^{-1}$ is trace class on $L^2 (-\infty, K)$.

Similarly, adding a boundary condition at $x=0$ is rank one, so with $H_D$ the
operator on $L^2 (0,K)$ with $u(0)=u(K)=0$ boundary conditions, it suffices to
prove that $f(H_D -z)^{-1}$ is trace class.

As $V\restriction [0,K]$ is in $L^1$, $H_D$ is bounded from below, and so by adding
a constant to $V$\!, we can suppose $H_D \geq 0$. Thus it suffices to show that
$f(H_D+1)^{-1}$ is trace class. Write
\begin{align*}
f(H_D+1)^{-1} &= \int_0^\infty e^{-t} fe^{-tH_D}\, dt \\
&= \int_0^\infty e^{-t} (fe^{-tH_D/2}) (e^{-tH_D/2})\, dt.
\end{align*}

By general principles (see \cite{Sxxi}), the integral kernel of $e^{-tH/2}$, call it
$P_t (x,y)$, obeys
\begin{alignat*}{2}
P_t (x,y) &\leq  Ct^{-1/2} \exp \biggl( \f{-(x-y)^2}{Dt}\biggr), \qquad && t\leq 1 \\
&\leq C, \qquad && t\geq 1.
\end{alignat*}
>From this it follows that for any $g\in L^2 (0,K)$,
\[
\|ge^{-sH_D}\|_2 \leq C\|g\|_{L^2} (1+\abs{s}^{-1/4})
\]
and $\|\cdot\|_2$ a Hilbert-Schmidt norm. Thus
\begin{align*}
\|f(H_D+1)^{-1}\|_1 &\leq \int_0^\infty e^{-t}\, \|fe^{-tH/2}\|_2 \,
\|e^{-tH/2}\|_2 \, dt \\
&< \infty.
\end{align*}

The proof of continuity and differentiability in $f$ follows from these estimates.
\end{proof}

\begin{remark} The use of Dirichlet decoupling and semigroup estimates to get trace
class results goes back to Deift--Simon \cite{DS}.
\end{remark}

\begin{corollary}\lb{C3.2} If $L_t$ is given by \eqref{1.30} and $z\in\bbC_+$, then
\begin{equation} \lb{3.3}
X_t =(L_t-z) (L_0-z)^{-1} -1\in \mathfrak{I}_1.
\end{equation}
Moreover, if $V(x)$ is continuous in a one-sided
neighborhood of $x=0$, $\tr (X_t)$ is differentiable at $t=0$ and
\begin{equation} \lb{3.4}
\left. \f{d}{dt}\, \tr  (X_t)\right|_{t=0} =V(0) G(0,0)
\end{equation}
where $G$ is the integral kernel of the operator $(L_0-z)^{-1}$.
\end{corollary}

\begin{proof} We have that
\begin{equation} \lb{3.5}
L_t -L_0 = V\chi_{[0,t]}
\end{equation}
so \eqref{3.3} follows from Proposition~\ref{P3.1}. By the continuity assumption,
$X_t$ has a piecewise continuous integral kernel $X_t (x,y) =V(x) \chi_{[0,t]}(x)
G(x,y)$, so (see, e.g., Theorem~3.9 in Simon \cite{STI})
\[
\tr (X_t) =\int_0^t V(x) G(x,x)\, dx
\]
from which \eqref{3.4} follows.
\end{proof}

The main result in this section is

\begin{theorem}\lb{T3.3} Let $V\in L_\loc^2 (0,\infty)$ with $\sigma_\ess (H)
=[0,\infty)$. Then for
\[
k\in\bbC_+ \backslash\{k=-i\kappa\mid -\kappa^2 \in \sigma (L_0)\}
\]
we have
\begin{equation} \lb{3.6}
a_x(k) =\det [(L_x -k^2)/(L_0 -k^2)].
\end{equation}
\end{theorem}

\begin{proof} By continuity, we can suppose $\Real k\neq 0$. Similarly, by
Proposition~\ref{P3.1}, we can suppose $V$ is continuous on $[0,x]$.

Let $\ti a_s(k)$ be the right-hand side of \eqref{3.6}. If we prove that for $0<t<x$,
\begin{equation} \lb{3.7}
\f{d}{dt} \, \log [\ti a_t(k)] = V(t) [ik + w(k;t)]^{-1}
\end{equation}
then, by \eqref{2.6} and $\ti a_0 (k)\equiv 1$, we could conclude \eqref{3.6}.

As
\[
\left. \f{d}{dA}\, \log\bigl[\det (1+A)\bigr]\right|_{A=0} =\tr (A)
\]
and for $t$ near $t_0$,
\begin{equation} \lb{3.8}
\log \ti a_t(k) =\log \ti a_{t_0}(k) + \log \det [(L_t -k^2)/(L_{t_0}-k^2)],
\end{equation}
\eqref{3.4} and \eqref{3.8} imply that
\[
\left. \f{d}{dt}\, \log [\ti a_t(k)]\right|_{t_0} =
V(t_0) G(t_0, t_0)
\]
where $G$ is the integral kernel for $(L_{t_0} -k^2)^{-1}$.  This leads to \eqref{3.7} after
writing the Green's function in terms of $\psi_+$ and $\psi_-^{(0)}$.
\end{proof}

\section{The Step-by-Step Sum Rule} \lb{s4}

In this section, we will prove a general step-by-step sum rule for all $V\in
L_\loc^2 ([0,\infty))$ that involves $\int_0^x V(y)^2 \, dy$. We begin with
a preliminary: Recall (see Proposition~\ref{P2.1}) that $\kappa_j(t)\geq 0$
is defined so that $\kappa_1(t)\geq\kappa_2(t)\geq \cdots$ and
$\{-\kappa_j(t)^2\}_{j=1}^{N(t)}$ are the negative eigenvalues of $L_t$
and $\kappa_j(t)=0$ if $j>N(t)$, which may be infinite.

\begin{proposition}\lb{P4.1} For any $V\in L_\loc^2 ([0,\infty))$ and $t\in (0,\infty)$,
\begin{equation} \lb{4.1}
\sum_j \, \abs{\kappa_j(t)^2 - \kappa_j(0)^2}<\infty .
\end{equation}
\end{proposition}

\begin{proof} Let $A(s) =-\f{d^2}{dx^2} + \chi_{[t,\infty)}V + s\chi_{[0,t)} V$
so $A(0)=L_t$ and $A(1) = L_0$. Let $E_j (s)$ denote the negative eigenvalues
$B(s)$ with $E_1\leq E_2 \leq\cdots$ and $E_j(s)=0$ if $j>N(s)$, the number of
negative eigenvalues of $A(s)$. Let $\psi_j(s)$ be the corresponding normalized
eigenvectors. Pick $a>0$ so for all $s\in [0,1]$, $A(s) \geq 1-a$.

By first-order eigenvalue perturbation \cite{RS4,Kato} (a.k.a.\ the Feynman--Hellman
theorem) if $j\leq N(s)$:
\begin{align*}
\f{d}{ds}\, E_j(s) &= \langle \psi_j(s), \chi_{[0,t)} V\psi_j(s)\rangle \\
&= (E_j(s) + a) \langle (A(s)+a)^{-1/2} \psi_j(s), \chi_{[0,t]}
V(A(s)+a)^{-1/2} \psi_j(s)\rangle
\end{align*}
so
\begin{equation} \lb{4.2}
\biggl| \f{dE_j(s)}{ds}\biggr| \leq 2a \abs{\langle\psi_j(s), (A(s)+a)^{-1/2}
\chi_{[0,t]} V(A(s)+a)^{-1/2} \psi_j(s)\rangle}
\end{equation}
and thus
\begin{align}
\sum_{j=1}^{N(s)} \, \biggl| \f{dE_j(s)}{ds}\biggr| &\leq 2a \|\chi_{[0,t]}
\abs{V}^{1/2}
(A(s)+a)^{-1/2}\|_2^2 \notag \\
&\leq C \lb{4.3}
\end{align}
where $\|\cdot\|_2$ is the Hilbert-Schmidt norm. \eqref{4.3}, which indicates a
$C$ independent of $s\in [0,1]$, follows from an estimate that
\begin{align*}
\|\chi_{[0,t]} \abs{V}^{1/2} (A(s)+a)^{-1/2}\|_2^2
&=\int_0^t \abs{V(x)}^2 (A(s)+a)^{-1} (x,x), dx \\
&\leq C.
\end{align*}
\eqref{4.3} implies
\[
\sum_{j=1}^\infty \, \abs{E_j(s) - E_j(u)} \leq C\abs{s-u}
\]
which for $s=1$, $u=0$ is \eqref{4.1}.
\end{proof}

\begin{remark}
One may also prove this proposition using the $\mathfrak{I}_1\to L^1$ bound for the Kre{\u\i}n spectral
shift function.  Indeed, the proof of this general result follows along the general lines given above.
\end{remark}

We can use this to define the Blaschke product that we will need to deal with the
zeros and poles of $a_t(k)$:

\begin{proposition} \lb{P4.2} Let
\begin{equation} \lb{4.4}
B_t(k) =\prod_j \biggl\{\biggl[ \f{k+i\kappa_j(0)}{k-i\kappa_j(0)}\,
\f{k-i\kappa_j(t)} {k+i\kappa_j(t)}\biggr] \exp \biggl[ -\f{2i}{k} \,
\bigl(\kappa_j(0)-\kappa_j(t)\bigr)\biggr] \biggr\}.
\end{equation}
Then:
\begin{SL}
\item The infinite product converges on $\bbC_+\backslash
\{\kappa_j^{(0)}\}_{j=0}^{N(0)}$.

\item $B_t(k)$ has a continuation to $\bar\bbC_+\backslash[\{\kappa_j(0)\}_{j=0}^N \cup \{0\}]$ and
\begin{equation} \lb{4.5}
k\in\bbR\backslash \{0\} \Rightarrow \abs{B_t(k)} =1.
\end{equation}
\item For $k\notin i\bbR$,
\begin{equation} \lb{4.6}
\abs{\log \abs{B_t(k)}} \leq C\abs{\Real k}^{-2}.
\end{equation}
\item Uniformly for $\arg(y)\leq \f{\pi}{4}$,
\begin{equation} \lb{4.6a}
\log [B_t (iy)]=\f{2}{3y^3}\, \sum_j \, [\kappa_j(0)^3 - \kappa_j(t)^3] + O(\abs{y}^{-5})
\end{equation}
as $\abs{y}\to\infty$.
\end{SL}
\end{proposition}

\begin{proof} Let $\kappa,\lambda >0$. Define
\begin{equation} \lb{4.7}
F(k;\kappa,\lambda) = \log \biggl[ \f{k+i\kappa}{k-i\kappa}\,
\f{k-i\lambda}{k+i\lambda}\biggr] -\f{2i}{k}\, (\kappa-\lambda).
\end{equation}
Then
\[
F(k; \lambda,\lambda)=0
\]
and, by a straightforward computation,
\begin{equation} \lb{4.8}
\f{\partial}{\partial\kappa}\, F(k; \kappa,\lambda) =
-\f{2i\kappa^2}{k(k^2 + \kappa^2)}.
\end{equation}
It follows for $k\in\bbC$ with $\pm ik\notin [\min(\kappa,\lambda)\max
(\kappa,\lambda)]$,
\begin{equation} \lb{4.9}
\abs{F(k;\kappa,\lambda)} \leq 2 \int_{\min (\kappa,\lambda)}^{\max(\kappa,\lambda)}
\f{\mu^2}{\abs{k}\, \abs{k^2 + \mu^2}}\, d\mu
\end{equation}
The right side is invariant under $k\to\bar k$, so suppose $\Ima k\geq 0$.
Then $\f{\mu}{\abs{k+i\mu}} \leq 1$, so
\begin{equation} \lb{4.10}
\abs{F(k;\kappa,\lambda)} \leq \f{\max(\kappa,\lambda)^2 -\min (\kappa,\lambda)^2}
{\abs{k}\, \inf \{\abs{k-i\mu}\mid \mu\in\pm (\min (\kappa,\lambda),
\max(\kappa, \lambda))\}}.
\end{equation}
We can thus prove:

\smallskip
(i) By \eqref{4.10}, if $k\notin \{0\}\cup \{i\kappa_j(0)\}\cup \{-i\kappa_j(t)\}
\equiv Q$, we have for all $n$ sufficiently large that
\[
\abs{F(k; \kappa_n(0), \kappa_n(t))} \leq C_k \abs{ \kappa_n(0)^2 - \kappa_n(t)^2}
\]
so, by \eqref{4.1}, the product \eqref{4.4} converges absolutely and uniformly
on compact subsets of $\bbC\backslash Q$.

\smallskip
(ii) The above argument shows $B$ has analytic continuation across $\bbR\backslash
\{0\}$. Since the continuation is given by a convergent product, and the finite
products have magnitude $1$ on $\bbR$,  that is true of $B$ on $\bbR\backslash \{0\}$.

\smallskip
(iii) From \eqref{4.10} and $\inf \{\abs{k-i\mu}\mid k\in \dots \}\geq \Real\abs{k}$,
we have
\[
\abs{F(k;\kappa,\lambda)} \leq \f{\abs{\kappa^2 -\lambda^2}}{\abs{\Real k}^2}
\]
which, given \eqref{4.1}, implies \eqref{4.6}.

\smallskip
(iv) By \eqref{4.8} for $y$ real and large,
\begin{align*}
\f{\partial}{\partial\kappa}\, F(iy, \kappa,\lambda) &= \f{2\kappa^2}{y(y^2 -\kappa^2)} \\
&= \f{2\kappa^2}{y^3} + O\biggl( \f{\kappa^4}{y^5}\biggr)
\end{align*}
so \eqref{4.6a} holds by integrating and using
\[
\int_{\kappa_j(t)}^{\kappa_j(0)} 2\mu^2 \, d\mu = \tfrac23\, [ \kappa_j(0)^3 - \kappa_j(t)^3 ].
\qedhere
\]
\end{proof}

Let $a_t (k)$ be given by \eqref{2.5} and $B_t (k)$ by \eqref{4.4}. The two functions
are analytic in $\bbC_+$ and have the same zeros and poles, so
\begin{equation} \lb{4.11}
g_t(k) = \log \biggl[ \f{a_t(k)}{B_t(k)}\biggr]
\end{equation}
is analytic in $\bbC_+$. We define $g_t$ by taking the branch of $\log$ which is real
for $k=i\kappa$ with $\kappa$ large.

\begin{proposition} \lb{P4.3}
\begin{SL}
\item $a_t (k)$ is analytic in $\bbC_+$.

\item For a.e.~$k\in\bbR_+$, $\lim_{\veps\downarrow 0} g_t (k+i\veps)
\equiv g_t (k)$ exists and if $\Im m(k^2+i0)>0$, then
\begin{equation} \lb{4.12}
\Real g_t (k) =\tfrac12\, \log \biggl[ \f{T(k,0)}{T(k,t)}\biggr]
\end{equation}
with $T$ given by \eqref{2.8}.

\item For each $\veps >0$,
\begin{equation} \lb{4.13}
\Ima k>\veps \Rightarrow \abs{g_t (k)}\leq C_\veps \abs{k}^{-1}.
\end{equation}

\item For all $k\in\bbC_+$, $\Real k\neq 0$,
\begin{equation} \lb{4.14}
\abs{g_t(k)}\leq C [\abs{\Real k}^{-1} + \abs{\Real k}^{-2}].
\end{equation}

\item As $y\to\infty$ along the real axis,
\begin{equation} \lb{4.15}
g_t (iy) =ay^{-1} + by^{-3} + o(y^{-3})
\end{equation}
with coefficients
\begin{align}
a &= -\tfrac12 \int_0^t V(x)\, dx \lb{4.16} \\
b &= \tfrac18 \int_0^t V(x)^2 \, dx - \tfrac23\, \sum_j \,
[\kappa_j(0)^3 - \kappa_j(t)^3].  \lb{4.17}
\end{align}
\end{SL}
\end{proposition}

\begin{proof} (i) is discussed in the definition.

\smallskip
(ii) This combines Proposition~\ref{P2.4} and \eqref{4.5}.

\smallskip
(iii) This follows from \eqref{2.6a}, \eqref{2.6b}, \eqref{4.6}, \eqref{4.6a}, and
the continuity (and so, boundedness) of $g_t$ on compact subsets of $\bbC_+$.

\smallskip
(iv) This combines \eqref{2.6a} and \eqref{4.6}.

\smallskip
(v) This combines \eqref{2.12} and \eqref{4.6a}.
\end{proof}

We are now ready for the nonlocal step-by-step sum rule.

\begin{theorem} \lb{T4.4} Suppose $V\in L_\loc^2 (\Reals^+)$ and $\Im m(E+i0)>0$ for almost every $E>0$.
Then for any $y_0, y_1 \in (0,\infty)$,
\begin{equation} \lb{4.18}
\Real \biggl[ g_t(iy_0) - \frac{y_1 g_t(iy_1)}{y_0}\biggr] =
     \int_0^\infty \f{(y_0^2 - y_1^2)\xi^2}{y_0(\xi^2 + y_0^2) (\xi^2 + y_1^2)}\,
        \log \biggl[ \f{T(\xi,0)} {T(\xi,t)}\biggr]\, \frac{d\xi}{\pi}
\end{equation}
where $g_t$ is given by \eqref{4.11} and $T$, by \eqref{2.8}.
\end{theorem}

\begin{proof} If $h$ is a bounded harmonic function on $\bbC_+$ with a continuous
extension to $\bar\bbC_+$, then for $y>0$,
\begin{equation} \lb{4.19}
h(x+iy) = \f{y}{\pi} \int \f{h(\xi)}{(\xi-x)^2 + y^2}\, d\xi .
\end{equation}
This Poisson representation is standard \cite{Rudin,Stein} and follows by noting that
the difference of the two sides is a harmonic function on $\bbC_+$ vanishing on
$\bbR$ so, by the reflection principle, a restriction of a bounded harmonic function
on $\bbC$ vanishing on $\bbR$ and so $0$ by Liouville's Theorem.

As $\Real g_t(k)$ is a bounded harmonic function on $\{k\mid\Ima k\geq \veps\}$,
we have for all $y>0$ and $\veps >0$,
\begin{equation} \lb{4.20}
\Real g_t (x+iy+i\veps) = \f{y}{\pi} \int \f{\Real g_t (\xi+i\veps)}{(\xi-x)^2 + y^2}\, dw
\end{equation}
and therefore,
\begin{equation} \lb{4.27}
\Re g_t(iy_0+i\veps) - \f{y_1}{y_0}\, \Re g_t (iy_1 + i\veps) =
    \int Q(\xi, y_0, y_1) \Real g_t (\xi+i\veps) \,d\xi
\end{equation}
where
\begin{align*}
Q(\xi, y_1, y_0) &= \f{1}{\pi} \biggl[ \f{y_0}{\xi^2 + y_0^2} - \f{y_1}{y_0}\, \f{y_1}{\xi^2 + y_1^2}\biggr] \\
&= \f{1}{\pi} \, \f{y_0^2 - y_1^2}{y_0} \, \f{\xi^2}{(\xi^2 + y_0^2)(\xi^2 + y_1^2)}.
\end{align*}

By \eqref{4.14}, uniformly in $\veps$,
\[
\abs{\Real g_t (\xi+i\veps)} \leq C [\abs{\xi}^{-2} + \abs{\xi}^{-1}]
\]
and clearly,
\[
\abs{Q(\xi)} \leq C_{y_0, y_1}\, \f{\xi^2}{1+\xi^4}
\]
so, by the dominated convergence theorem, we can take $\veps\downarrow 0$ in
\eqref{4.20}. The left side converges to the left side of \eqref{4.18} and, by
\eqref{4.12} and $\Real g_t (-\bar k)=\Real g_t (k)$, the right side converges to
the right side of \eqref{4.18}.
\end{proof}

Here is the step-by-step version of the Faddeev--Zhabat sum rule \eqref{1.27}:

\begin{theorem}[Step-by-Step Faddeev--Zhabat Sum Rule]\lb{T4.5} Suppose $V\in L_\loc^2(\Reals^+)$ and
$\Im m(E+i0)>0$ for almost every $E>0$. For any $t>0$,
\begin{equation} \lb{4.28}
\tfrac18 \int_0^t V(x)^2\, dx = \tfrac23 \sum_j\, [\kappa_j(0)^3 - \kappa_j(t)^3]
    + \lim_{y\to\infty} \int_0^\infty P(\xi,y)
\log \biggl[ \f{T(\xi,t)}{T(\xi,0)}\biggr] \, d\xi
\end{equation}
where
\begin{equation} \lb{4.29}
P(\xi,y) = \f{1}{\pi} \biggl[ \f{4\xi^2 y^4}{(\xi^2 + 4y^2)(\xi^2 +y^2)}\biggr].
\end{equation}
\end{theorem}

\begin{proof} By \eqref{4.15}, with $b$ given by \eqref{4.17},
\[
y^3 [g_t (2iy) -\tfrac12\, g_t(iy)] = b[\tfrac18 -\tfrac12] + o(1)
\]
so, by \eqref{4.18},
\begin{equation} \lb{4.30}
b=\lim_{y\to\infty} - \f83 \, \biggl[ \f{y^3 [(2y)^2 -(y)^2]}{2\pi y}\biggr]
\int_0^\infty \f{\xi^2}{(\xi^2+4y^2)(\xi^2 + y^2)} \, \log\biggl[\f{T(\xi,0)}{T(\xi,t)}
\biggr]\, d\xi
\end{equation}
which is \eqref{4.28}.
\end{proof}

\begin{remarks} 1. As $\lim_{y\to\infty} P(\xi,y) =\f{1}{\pi}\xi^2$, formally, \eqref{4.28} is
just a difference of \eqref{1.28} for $L_0$ and $L_t$.

2. In the preceding theorems, the assumption that $\Im m(E+i0)>0$ for almost every $E>0$ was only used
to allow us to apply Proposition~\ref{P2.4} to obtain a simpler expression for the boundary values of $a_t(k)$.
The assumption may be removed if one is willing to replace the ratio $T(\xi,t)/T(\xi,0)$ by the limiting value
of the relative Wronskian.
\end{remarks}

\section{Lower Semicontinuity of the Quasi-Szeg\H{o} Terms} \lb{s5}

For any $V\in L_\loc^1 (0,\infty)$, we can define (in the limit circle case after
picking a boundary condition at infinity) $T(k,0)$ by \eqref{2.8} for a.e.\
$k\in (0,\infty)$ and then
\begin{equation} \lb{5.1}
Q(V) =-\f{1}{\pi} \int_0^\infty \log [T(k,0)] k^2 \, dk.
\end{equation}
Since $T\leq 1$, $-\log [T]\geq 0$ and the integral can only diverge to $\infty$,
so $Q(M)$ is always defined although it may be infinite. The main result in this
section is:

\begin{theorem}\lb{T5.1} Let $V_n, V$ be a sequence in $L_\loc^2 ((0,\infty))$.
Let $V$ be limit point at infinity. Suppose
\begin{equation} \lb{5.2}
\int_0^a \abs{V_n(x)-V(x)}^2\, dx\to 0
\end{equation}
for each $a>0$. Then
\begin{equation} \lb{5.3}
Q(V) \leq \liminf Q(V_n).
\end{equation}
\end{theorem}

\begin{remarks} 1. As noted in the Introduction, this is related to results in
Sylvester--Winebrenner \cite{SW}. However, they have no bound states and $\abs{r(k)}
\leq 1$ in the upper half-plane. This fails in our case and our argument will need
to be more involved.

2. It is interesting that the analogue in the Jacobi case \cite{KS} used
semicontinuity of the entropy and this result comes from weak semicontinuity of
the $L^p$-norm.

3. It is not hard to see that this result holds if $L_\loc^2$ is replaced by
$L_\loc^1$ and the $\abs{\dots}^2$ in \eqref{5.2} is replaced by $\abs{\dots}^1$.
Basically, one still has strong resolvent convergence in that case. But the argument
is simpler in the $L_\loc^2$ case we need, so that is what we state.
\end{remarks}

We will prove this theorem in several steps.  We will write $w_n(k)$ and $w(k)$ for the
$m$-functions (parameterized by momentum) associated to $V_n$ and $V$ respectively.

\begin{proposition}\lb{P5.2} Let $V_n,V$ obey the hypothesis of Theorem~\ref{T5.1}.
Then for all $k$ with $\Real k>0$ and $\Ima k>0$, one has $w_n(k)\to w(k)$ as $n\to\infty$.
\end{proposition}

\begin{proof} Let $H$ (resp., $H_n$) be the operator $u\mapsto -u'' +Vu$ on
$L^2 (0,\infty)$ with boundary condition $u(0)=0$ at $x=0$ and, if need be,
a boundary condition at $\infty$ for some $n$ if the corresponding $H_n$ is
limit circle at $\infty$.

By the standard construction of these operators, $H$ being limit point at
infinity has $D\equiv \{u\in C_0^\infty ([0,\infty))\mid u(0)=0\}$ as an
operator core. (\cite[Theorem~X.7]{RS2} has the result essentially if $V$ is
continuous, but the proof works if $V$ is $L_\loc^2$. Essentially, any $\varphi
\in [(H+i) [D]]^\perp$ solves $-\varphi''+ V\varphi =i\varphi$ with $\varphi(0)
=0$ and that cannot be $L^2$; it follows that $\ol{H\pm i[D]}=L^2$ which is
essential selfadjointness.)

Let $f=(H-k^2)\varphi$ with $\varphi\in D$. Then
\begin{align*}
\|[(H_n -k^2)^{-1} -(H-k^2)^{-1}] f\| &= \|(H_n-k^2)^{-1} (V_n-V)\varphi\| \\
&\leq \abs{\Ima k^2}^{-1} \|(V_n-V)\varphi\|\to 0
\end{align*}
by \eqref{5.2}, so we have strong resolvent convergence.

If $\varphi\in L^2 (0,a)$ and $\psi =(H_n-k^2)^{-1}\varphi$, then for $x>a$,
\[
w_n (k,x) = \f{\psi'_n(x)}{\psi_n(x)}
\]
and so, for $x>0$, we have $w_n(k,x)\to w(k,x)$.

Differentiating \eqref{E:w} with respect to $x$ and using \eqref{R2.1} leads to the Riccati equation
\begin{equation}\label{E:Riccati}
\frac{d w}{d x} = k^2 - V(x) - w^2.
\end{equation}
By combining this with \eqref{5.2}, one can deduce $w_n(k)\to w(k)$.
\end{proof}

We now define the reflection coefficient (for now, a definition; we will
discuss its connection with reflection at the end of the section) by
\begin{equation} \lb{5.4}
r_n (k) = \f{ik-w_n(k)}{ik+w_n(k)}.
\end{equation}
The following bound is clearly relevant.

\begin{proposition}\lb{P5.3} Let $k=\abs{k}e^{i\eta}$ with $\eta\in
[0,\f{\pi}{2})$, $\abs{k}\neq 0$. Then
\begin{equation} \lb{5.5}
\sup_{z\in\bbC_+}\, \biggl| \f{ik-z}{ik+z}\biggr| =
\biggl( \f{1+\sin (\eta)}{1-\sin(\eta)}\biggr)^{1/2}.
\end{equation}
\end{proposition}

\begin{proof} $z\mapsto\f{ik-z}{ik+z}$ is a fractional linear transformation
which takes $z = -ik\in \bbC_-$ to infinity since $\Real k>0$ if $\eta\in
[0,\f{\pi}{2})$. Thus $\bbC_+$ is mapped into the interior of the circle
$\{\f{ik-x}{ik+x}\mid x\in\bbR\}\cup\{-1\}$. By replacing $k$ by $k/\abs{k}$,
we can suppose $\abs{k}=1$. Let
\[
f(x) =\biggl| \f{x-ie^{i\eta}}{x+ie^{i\eta}}\biggr|^2.
\]

Straightforward calculus shows that $f'(x)=0$ exactly at $x=\pm 1$. Since
$\abs{f}\to 1$ as $x\to \pm\infty$, we see the maximum of $f(x) =(1+x^2
+ 2x\sin\eta)/(1+x^2-2x\sin\eta)$ occurs at $x=1$ and is $(1 +
\sin(\eta))/(1-\sin(\eta))$.
\end{proof}

\begin{lemma} \lb{L5.4} Let $f_n$ and $f_\infty$ be a sequence of functions on
$\bbD$, the open disk, with \begin{equation} \lb{5.6}
\sup_{z\in\bbD,n}\, \abs{f_n(z)}<\infty.
\end{equation}
Let $f_n(z)\to f_\infty(z)$ for all $z\in\bbD$. Let $f_n (e^{i\theta})$ be the
a.e.\ radial limit of $f_n (re^{i\theta})$ and similarly for $f_\infty (e^{i\theta})$.
Then $f_n (e^{i\theta})\to f_\infty (e^{i\theta})$ weak-$*$, that is, for all $g\in
L^1 (\partial\bbD)$,
\begin{equation} \lb{5.7}
\int_0^{2\pi} g(e^{i\theta}) f_n (e^{i\theta}) \, \f{d\theta}{2\pi} \to
\int_0^{2\pi} g(e^{i\theta})
f_\infty (e^{i\theta}) \, \f{d\theta}{2\pi}.
\end{equation}
\end{lemma}

\begin{proof} By \eqref{5.6}, it suffices to prove \eqref{5.7} for $g(e^{i\theta})
= e^{ik\theta}$ for all $k$. But for $H^\infty $ functions (see \cite{Rudin}),
$\int e^{ik\theta} f(e^{i\theta}) \f{d\theta}{2\pi} =0$ if $k>0$ and $\int
e^{-ik\theta} f(e^{i\theta}) \f{d\theta}{2\pi} = f^{(k)}(0)/k!$. Pointwise
convergence in $\bbD$ and boundedness implies convergence of all derivatives
inside $\bbD$.
\end{proof}

\begin{theorem}\lb{T5.5} Let $r_n(k)$ be given by \eqref{5.4} for $\Ima k>0$.
Then for a.e.~$k\in (0,\infty)$, $r_n(k)=\lim_{\veps\downarrow 0} r_n (k+i\veps)$
exists and obeys
\begin{equation} \lb{5.8}
\abs{r_n(k)} \leq 1, \qquad (k>0).
\end{equation}
Moreover, for any $g$ in $L^1 (a,b)$ with $0<a<b<\infty$, we have that
\begin{equation} \lb{5.9}
\int_a^b g(k) r_n(k)\, dk \to \int_a^b g(k) r(k)\, dk
\end{equation}
and that for $1\leq p<\infty$,
\begin{equation} \lb{5.10}
\liminf_{n\to\infty}\, \int_a^b \abs{r_n(k)}^p k^2 \, dk \geq \int_a^b
\abs{r(k)}^p k^2 \, dk.
\end{equation}
\end{theorem}

\begin{proof} Pick $0<c<a<b<d<\infty$. Let $Q$ be the semidisk in $\bbC_+$
with flat edge $(c,d)$. Let $\varphi:\bbD\to Q$ be a conformal map. Since
\[
\sup_{k\in Q}\, \arg(k) < \f{\pi}{2},
\]
we have
\begin{equation} \lb{5.11}
\sup_{n,k\in Q}\, \abs{r_n(k)}<\infty
\end{equation}
by Proposition~\ref{P5.3}. We can thus apply Lemma~\ref{L5.4} to $r_n \circ
\varphi$ and so conclude \eqref{5.9}. \eqref{5.8} follows from
Proposition~\ref{P5.3} for $\eta=0$.

Note that \eqref{5.9} implies $r_n\to r$ in the weak topology on $L^p ((a,b),k^2\, dk)$.
Thus \eqref{5.10} is just an expression of the fact that the norm
on a Banach space is weakly lower semicontinuous.
\end{proof}

\begin{proof}[Proof of Theorem~\ref{T5.1}] Notice that
\begin{equation} \lb{5.12}
T(k_0,0)+\abs{r(k_0,0)}^2 =1.
\end{equation}
Thus
\begin{align}
-\log [T] &=-\log (1-\abs{r}^2) \notag \\
&= \sum_{m=1}^\infty \, \f{\abs{r}^{2m}}{m}. \lb{5.13}
\end{align}
\eqref{5.10} implies that for each $m$ and $0<a<b<\infty$,
\[
\int_a^b \biggl[\f{\abs{r(k)}^{2m}}{m}\biggr] k^2 \, dk \leq
\liminf \int_a^b \biggl[ \f{\abs{r_n(k)}^{2m}}{m}\biggr] k^2\, dk,
\]
which becomes
\[
\int_a^b \biggl[\,\sum_{m=1}^M \, \f{\abs{r}^{2m}}{m}\biggr]
k^2\, dk \leq \liminf \int_a^b \biggl[\, \sum_{m=1}^M \,
\f{\abs{r_n}^{2m}}{m}\biggr] k^2\, dk
\]
so, by \eqref{5.13},
\[
-\f{1}{\pi} \int_a^b \log [T(k_0,0)]k^2, dk \leq
\liminf \biggl\{-\f{1}{\pi} \int_0^\infty \log
[T_n (k_0,0)] k^2\, dk\biggr\}.
\]
Now take $a\downarrow 0$ and $b\to\infty$.
\end{proof}

\smallskip
We end this section with a sketch of an alternate approach to
Theorem~\ref{T5.1}. We present this approach because it is rooted in
the physics of scattering. Since we have a direct proof, we do not
produce all the technical details---indeed, one is missing. The
argument is in a sequence of steps:

\smallskip
\noindent{\bf Step 1.} \ Let $L$ be the whole-line problem obtained by
setting $V=0$ on $(-\infty, 0)$. Let $j$ be a $C^\infty$ function with
$0\leq j\leq 1$ and $j(x)=0$ if $x>0$ and $j(x)=1$ if $x<-1$. Let $J$
be multiplication by $j$. Then, by \cite{Sim99},
\begin{equation} \lb{5.14}
\slim_{t\to\pm\infty}\, e^{itL} Je^{-itL} P_\ac (L) = P_\ell^\pm (L)
\end{equation}
exist and are invariant projections for $L$. $L\restriction\ran (P_\ell^\pm)$
is absolutely continuous and has spectrum $[0,\infty)$ with multiplicity $1$.

\smallskip
\noindent{\bf Step 2.}
\begin{equation} \lb{5.15}
P_\ell^- (L) P_\ell^+(L) P_\ell^- (L) \equiv R_\ell^- (L)
\end{equation}
is a positive operator on $\ran (P_\ell^-)$ which commutes with $L\restriction
\ran (P_\ell^- (L))$ and so, by the simplicity of the spectrum of this operator,
it is multiplication by a function $R_L(E)$. Since $0\leq R_\ell^-(L)\leq 1$,
as a function, $0\leq R(E)\leq 1$. $R$ is discussed in \cite{Sim99}.

\smallskip
\noindent{\bf Step 3.} \ By computations related to those in \cite{SW},
\begin{equation} \lb{5.16}
R_L (k) = \abs{r(k)}^2
\end{equation}
with $r$ given by \eqref{5.4}.

\smallskip
\noindent{\bf Step 4.} \ We believe that for $V_n\to V$ in the sense of
Theorem~\ref{T5.1}, one has for a dense set of vectors uniformity in $n$ of
the limit in \eqref{5.14}, but we have not nailed down the details. If true,
one has
\begin{equation} \lb{5.17}
\wlim_{n\to\infty} \, R_\ell^-(L_n) = R_\ell^- (L).
\end{equation}

\smallskip
\noindent{\bf Step 5.} \ By \eqref{5.16}, $\abs{r_n(k)}^2 \to \abs{r(k)}^2$ weakly
as $L^\infty$-functions (i.e., when smeared with $g\in L^1 (a,b)$) on $[a,b]$ for
any $0<a<b<\infty$. By the weak semicontinuity of the norm, \eqref{5.10} holds for
$p\geq 2$.

\smallskip
\noindent{\bf Step 6.} \ Get semicontinuity of $Q(V)$ from \eqref{5.10} for
$p\geq 2$, as we do in the above proof.

\section{Local Solubility} \lb{s6}

In this section, we will study \eqref{1.19x} and describe its relation to $d\rho$ being
the spectral measure of some $V\in L_\loc^2$. We will prove:

\begin{theorem}\lb{T6.1} Let $d\rho$ be a measure obeying condition
{\rm{(i)}} of Theorem~\ref{T1.3}. Define $F$ by \eqref{1.18} and suppose
\eqref{1.19x} holds. Then $d\rho$ is the spectral measure of some $V\in L_\loc^2$.
\end{theorem}

\begin{theorem}\lb{T6.2} Let $d\rho$ be the spectral measure of a potential in $L^2$.
Then \eqref{1.19x} holds, that is, $F\in L^2(\Reals^+)$.
\end{theorem}

Before discussing the main ideas used to prove these results, we wish to reassure the reader
that the hypotheses of Theorem~\ref{T6.1} do bound the growth of $d\rho$ at infinity.  Specifically,
we know that \eqref{HRTa} must hold for any spectral measure.  We do this first because such information
is helpful in justifying some calculations that appear once the real work begins.

\begin{lemma}
If $d\rho$ obeys condition {\rm{(i)}} of Theorem~\ref{T1.3} and \eqref{1.19x} holds, then
\begin{equation}\label{HRTaAgain}
\int \frac{d\rho(E)}{1+E^2} < \infty.
\end{equation}
\end{lemma}

\begin{proof}
Unravelling the definitions of $F(q)$ and $d\nu$ given in \eqref{1.18} and \eqref{1.17a}, we find
$$
F(q) = \pi^{-1/2} \int_1^\infty \exp\bigl\{ -\bigl(q-\sqrt{E}\bigr)^2 \bigr\} E^{-1}\,d[\rho-\rho_0](E).
$$
The contribution of $\rho_0$ can be bounded using
$$
\tfrac{1}{\pi} \int_0^\infty \exp\bigl\{ -\bigl(q-\sqrt{E}\bigr)^2 \bigr\} E^{-1/2} \,dE
= \tfrac{2}{\pi} \int_0^\infty \exp\bigl\{ -(q-k)^2 \bigr\}\,dk \leq 2 \pi^{-1/2},
$$
which shows that
$$
\int_1^\infty \exp\bigl\{ -\bigl(q-\sqrt{E}\bigr)^2 \bigr\} E^{-1} \,d\rho(E) \leq 2 + 2|F(q)|.
$$
Integrating both sides $\frac{dq}{1+q^2}$ leads to \eqref{HRTaAgain}, at least when the region of
integration is restricted to $[1,\infty)$.  The remaining portion of the integral is finite by condition
{\rm{(i)}} of Theorem~\ref{T1.3}.
\end{proof}

The key to proving the two theorems of this section will be the fact that essentially, $\widehat
F(\alpha)$, the Fourier transform of $F$, is $e^{-\f14 \alpha^2} A(\alpha)$, where $A(\alpha)$ is
the $A$-function introduced by Simon \cite{Sim271} and studied further by Gesztesy--Simon
\cite{GSAnn}.

We will, first and foremost, use formula (1.21) from \cite{GSAnn}:
\begin{equation} \lb{6.1}
A(\alpha) =-2\int \lambda^{-1/2} \sin \bigl(2\alpha \sqrt{\lambda}\bigr) \, d(\rho-\rho_0)
(\lambda)
\end{equation}
where $\lambda^{-1/2}\sin (2\alpha \sqrt{\lambda})$ is interpreted as $\abs{\lambda}^{-1/2}
\sinh (2\alpha \sqrt{\abs{\lambda}})$ if $\lambda <0$ and \eqref{6.1} holds in
distributional sense. We will also need the following (eqn.~(1.16) of \cite{GSAnn}):
\begin{equation} \lb{6.2}
\abs{A(\alpha)-V(\alpha)} \leq \biggl| \int_0^\alpha \abs{V(y)}\,dy\biggr|^2
\exp \biggl( \alpha \int_0^\alpha \abs{V(y)}\, dy \biggr)
\end{equation}
proven in \cite{Sim271} for regular $V$'s and in (1.16) of \cite{GSAnn} for
$V\in L_\loc^1$. Finally, we need the following result, which follows readily from Remling's work \cite{Rem02,Rem03}.
(It can also be proved using the Gel'fand--Levitan method.)

\begin{proposition}\lb{P6.3} Let $d\rho$ be a measure obeying \eqref{HRTaAgain} and condition {\rm{(i)}}
of Theorem~\ref{T1.3}. If the distribution \eqref{6.1} lies in $L_\loc^1 [0,\infty)$, then $d\rho$ is the spectral
measure of a potential $V\in L_\loc^1 [0,\infty)$.
\end{proposition}

\begin{proof}
Consider the continuous function
$$
K(x,t)=\tfrac12\phi(x-t)-\tfrac12\phi(x+t),
\quad\text{where}\quad
\phi(x)=\int_0^{|x|/2} A(\alpha)\,d\alpha
$$
and $A(\alpha)$ is given by \eqref{6.1}.
As explained in Theorem~1.1 of \cite{Rem03}, $A(\alpha)$ is the $A$-function of a potential in $L^1_\loc$ provided
\begin{equation} \label{PosDef}
\int\int \bar\psi(x)\psi(t) \bigl[\delta(x-t) + K(x,t)\bigr]\,dx\,dt > 0
\end{equation}
for all non-zero $\psi\in L^2([0,\infty))$ of compact support.  We will now show that this holds.

For $\psi\in C_c^\infty$, elementary manipulations using \eqref{6.1} show
\begin{equation*}
\int\!\!\int \bar\psi(x)\psi(t) K(x,t) \,dx\,dt
= \int\!\!\int\!\!\int \bar\psi(x)\psi(t) \frac{\sin(x\sqrt{\lambda})\sin(t\sqrt{\lambda})}{\lambda}
    \,dx\,dt\,d[\rho-\rho_0](\lambda).
\end{equation*}
Thus by recognizing the spectral resolution of the free Schr\"odinger operator we have
\begin{equation*}
\text{LHS\eqref{PosDef}} = \int \biggl| \int \psi(x) \frac{\sin(x\sqrt{\lambda})}{\sqrt{\lambda}}\,dx \biggr|^2
     \,d\rho(\lambda)
\end{equation*}
for such test functions.  It then extends easily to all $\psi\in L^2([0,\infty))$ of compact support, because $K$
is a bounded function.

This representation shows that LHS\eqref{PosDef} is non-negative.  It cannot vanish for non-zero
$\psi$ because the Fourier sine transform of $\psi$ is analytic and so has discrete zeros; however,
the support of $d\rho$ is not discrete by hypothesis. Thus we have shown that $A(\alpha)$ defined
by \eqref{6.1} is the $A$-function of some $V\in L^1_\loc$.

Unfortunately, we are only half-way through the proof; the $A$-function need not uniquely determine the spectral
measure through \eqref{6.1}.  This is the case, for example, when the potential is limit circle at infinity; different
boundary conditions lead to different spectral measures, but all have the same $A$-function.  Christian Remling has
explained to us that using de Branges work, \cite{deB}, one can deduce that this is actually the only way
non-uniqueness can occur.  In our situation however, we have some extra information which permits us to complete
the proof of uniqueness without much technology, which is what we proceed to do now.

Let $d\rho_1$ denote the spectral measure for the potential $V$ just constructed (with a boundary condition at infinity
if necessary).  Classical results tell us that \eqref{HRTaAgain} holds for $d\rho_1$ and that
$\int_{-\infty}^0 \exp\{c\sqrt{-\lambda}\}\,d\rho_1(\lambda)<\infty$ for any $c>0$.  Lastly, by construction
we have
\begin{equation} \label{Aequal}
\int \lambda^{-1/2} \sin \bigl(2\alpha \sqrt{\lambda}\bigr) \, d(\rho-\rho_0)(\lambda)
=\int \lambda^{-1/2} \sin \bigl(2\alpha \sqrt{\lambda}\bigr) \, d(\rho_1-\rho_0)(\lambda)
\end{equation}
as weak integrals of distributions.  We wish to conclude that $\rho_1=\rho$.

Our first step is to prove that the support of $d\rho_1$ is bounded from below.  Let us fix
a non-negative $\phi\in C^\infty_c(\Reals)$ with $\int \phi(x)\,dx=1$ and $\supp(\phi)\subset[1,2]$.
Elementary considerations show that there is a constant $C$ so that
$$
\biggl| \int k^{-1} \sin(2\alpha k) \phi(\alpha/N)\,d\alpha \biggr| \leq C N^2 (1+k)^{-100}
$$
for all $N>1$ and all $k\geq 0$.  More easily, we have
$$
4N^2 e^{4Nk} \geq \tfrac{N}{k} \sinh(4Nk) \geq \int k^{-1} \sinh(2\alpha k) \phi(\tfrac\alpha{N})\,d\alpha
\geq \tfrac{N}{k} \sinh(2Nk) \geq N^2 e^{Nk}
$$
for the same range of $N$ and $k$.  Putting this together with \eqref{Aequal} we obtain
\begin{equation}\label{rho1bb}
\int_{-\infty}^0  e^{N\sqrt{-\lambda}} \, d\rho_1(\lambda) \leq C_1 + C_2 e^{4N\sqrt{-E_1}},
\end{equation}
where $E_1$ denotes the infimum of the support of $d\rho$ just as in condition {\rm{(i)}}
of Theorem~\ref{T1.3}.  Taking $N\to\infty$ in \eqref{rho1bb} leads to the conclusion that
the support of $\rho_1$ is bounded from below (by $16E_1$, which is easily improved).

Now that we know that the supports of both $\rho$ and $\rho_1$ are bounded from below, we may use
\begin{equation} \label{FouInt}
\tfrac{2}{\sqrt{\pi}}\int \alpha e^{-\alpha^2/s} \sin(2\alpha k)\,d\alpha = s^{3/2} k\, e^{-s k^2},
    \quad\text{for $s>0$ and $k\in\Cmplx$},
\end{equation}
on both sides of \eqref{Aequal} and so obtain
\begin{equation} \label{LaplaceEqual}
\int e^{-s\lambda} \, d\rho(\lambda) = \int e^{-s\lambda} \, d\rho_1(\lambda).
\end{equation}
That $\rho_1=\rho$ now follows from the invertibility of Laplace transforms.
\end{proof}

As outlined above, our discussion of the local solubility condition revolves around a relation between the
distributions $A$ and $F$. Let
\begin{equation} \lb{6.3}
A(\alpha) = A_S(\alpha) + A_L (\alpha)
\end{equation}
where $A_S$ is the integral over $\lambda <1$ and $A_L$ over $\lambda\geq 1$. Since
\begin{equation} \lb{6.4}
\int_{-\infty}^\infty \pi^{-1/2} e^{-q^2} e^{-iq\alpha} = e^{-\alpha^2/4},
\end{equation}
\eqref{1.17a}, \eqref{1.18}, and \eqref{6.1} immediately imply
\begin{equation} \lb{6.5}
e^{-\alpha^2/4} A_L(\alpha) = i [\widehat F(2\alpha) - \widehat F(-2\alpha)].
\end{equation}

For $p\geq 1$ and $q\leq 0$, we have $e^{-(p-q)^2} \leq e^{-p^2} e^{-q^2}$.
Combining this with
\[
\int_{p\geq 1} e^{-p^2} d\abs{\nu} (p) <\infty,
\]
which follows from \eqref{HRTaAgain}, we obtain that for $q\leq 0$,
\begin{equation} \lb{6.6}
 F(q) \leq Ce^{-q^2}.
\end{equation}

\begin{proof}[Proof of Theorem~\ref{T6.1}] By \eqref{1.19x} and \eqref{6.6}, $F\in L^2(\bbR)$
and hence $\widehat{F}\in L^2(\bbR)$. By \eqref{6.5}, $A_L(\alpha)\in L_\loc^2$. By
\eqref{6.1}, $A_S(\alpha)$ is bounded on bounded intervals, so $A(\alpha)\in L_\loc^2$.
By Remling's Theorem (Proposition~\ref{P6.3}), $d\rho$ is the spectral measure of some
$V\in L_\loc^1$. By \eqref{6.2}, $\abs{A(\alpha)-V(\alpha)}$ is bounded on bounded
intervals, so $A\in L_\loc^2 \Rightarrow V\in L_\loc^2$.
\end{proof}

To prove Theorem~\ref{T6.2}, we need the following elementary fact:

\begin{proposition}\lb{P6.6} If $T$ is a tempered distribution on $(1,\infty)$ which is real
and $\Ima \widehat T(\alpha)\in L^2$, then $T\in L^2$.
\end{proposition}

\begin{proof}
We begin by noting that if $h\in L^2 (0,\infty)$, then
\begin{equation} \lb{6.7}
\int_{-\infty}^\infty \abs{\Real \widehat h(\alpha)}^2\,d\alpha =
\int_{-\infty}^\infty \abs{\Ima \widehat h(\alpha)}^2\,d\alpha
\end{equation}
since $\ol{\widehat h(\alpha)} =\widehat{\ti h}(\alpha)$, where $\ti h(x)=h(-x)$, and
thus, by the Plancherel theorem,
\begin{equation} \lb{6.8}
\int_{-\infty}^\infty h(\alpha)^2\,d\alpha = \int_{-\infty}^\infty \ti h(x)h(x)\,dx =0
\end{equation}
implying \eqref{6.7}. In particular,
\begin{equation} \lb{6.9}
\int_{-\infty}^\infty \abs{\widehat h(\alpha)}^2\, d\alpha = 2\int
\abs{\Ima \widehat h(\alpha)}^2\, d\alpha.
\end{equation}

Given $T$, pick $C^\infty$ $g$ with $g(x)=g(-x)$, $\abs{g(x)}\leq 1$, $g(x)=1$ for
$\abs{x}$ small, $\supp(y)\subset [-1,1]$, and $\int g(x)\,dx =1$. Define
\begin{equation} \lb{6.10}
g_L(x)=g \biggl( \f{x}{L}\biggr) \qquad
r_\delta(x)=\delta^{-1} g\biggl(\f{x}{\delta}\biggr)
\end{equation}
and note that $r_\delta * (g_LT)\in L^2$, supported in $(0,\infty)$ for $\delta <1$,
and since
\begin{equation} \lb{6.11}
[r_\delta * (g_L T)]\widehat{\ } (\alpha) = \widehat h_1(\alpha\delta) (\widehat g_L *T)(\alpha)
\end{equation}
and $\widehat h_1, \widehat g_L$ are real, we have
\[
\int \abs{\Ima [r_\delta * (g_L T)]^\sim(\alpha)}^2\, d\alpha \leq
\int \abs{\Ima \widehat T}^2\, d\alpha.
\]
Thus, by \eqref{6.9} and the Plancherel theorem,
\[
\int \abs{r_\delta * g_L T(x)}^2\, dx \leq 2\int \abs{\Ima \widehat T(\alpha)}^2 \,
d\alpha
\]
so $T\in L^2$ by taking $\delta\downarrow 0$ and $L\to\infty$.
\end{proof}

\begin{proof}[Proof of Theorem~\ref{T6.2}] If $V\in L^2$, then
\begin{equation} \lb{6.12}
\int_0^\alpha \abs{V(y)}\, dy \leq \biggl( \int_0^\infty \abs{V(y)}^2\, dy \biggr)^{1/2}
\alpha^{1/2}
\end{equation}
so \eqref{6.2} says that
\begin{equation} \lb{6.13}
\abs{A(\alpha) - V(\alpha)} \leq C \alpha^2 \exp (C\alpha^{3/2})
\end{equation}
and thus, $e^{-\alpha^2/2} A(\alpha)\in L^2$.

{}From \eqref{6.1},
\begin{equation} \lb{6.14}
\abs{A_S(\alpha)}\leq e^{C\alpha}
\end{equation}
so $e^{-\alpha^2/2} A_S(\alpha)\in L^2$, and thus, $e^{-\alpha^2/2} A(\alpha)\in L^2$.
By \eqref{6.5} and the fact that $F$ is real-valued, it follows that $\Ima \widehat{F}\in L^2$.

$F$ is not supported on $(1,\infty)$, but by \eqref{6.6} and boundedness on $(0,1)$,
$F=F_1 + F_2$, where $F_2$ is supported on $(1,\infty)$ and $F_1\in L^2$. Thus,
$\Ima \widehat{F}_1\in L^2$, so $\Ima\widehat{F}_2\in L^2$. By
Proposition~\ref{P6.6}, $F_2\in L^2$, that is, \eqref{1.19x} holds.
\end{proof}

\section{Proof of Theorem~\ref{T1.3}} \lb{s7}

Here we will use the results of the last three sections to prove Theorem~\ref{T1.3}. We
use the strategy of \cite{KS} as refined in \cite{SZ} and \cite{Sim288}.  We treat each
direction of the theorem in a separate subsection.

\subsection*{$\bf V\in L^2 \Rightarrow$ (i)--(iv)}  As $V\in L^2$, $V(H_0+1)^{-1}$
is compact, and thus (i) holds by Weyl's Theorem. (ii) is just Theorem~\ref{T6.2}.

Fix $R<\infty$ and let
\begin{equation} \lb{7.1}
V^{(R)}(x) = \begin{cases} V(x), & 0 \leq x \leq R \\
0, & x>R \end{cases}
\end{equation}
so $L_t=H_0$ if $t>R$. Thus applying Theorem~\ref{T4.5} to $V^{(R)}$ with $t>R$ gives
\begin{equation}\lb{7.2}
\begin{aligned}
\tfrac12 \int_0^R V(x)^2  &= \tfrac23 \sum_j [\kappa_j^{(R)}]^3 + \lim_{y\to\infty} \int_0^\infty
P(\xi,y)\log \biggl( \f{1}{T^{(R)}(\xi,0)}\biggr)\,d\xi.
\end{aligned}
\end{equation}
By rewriting $P$ as
\[
P(\xi,y) = \f{\xi^2}{\pi} \biggl[\f{1}{(1+\f{\xi^2}{4y^2}) (1+\f{\xi^2}{y^2})}\biggr],
\]
we see that it is monotone increasing in $y$. As the integrand $\log (\f{1}{T}) \geq 0$, the monotone
convergence theorem implies
\begin{equation} \lb{7.3}
\tfrac18 \int_0^R V(x)^2\, dx = \tfrac23\, \sum_{j=1}^{N(R)} \bigl[\kappa_j^{(R)}\bigr]^3 + Q(V^{(R)}).
\end{equation}

Now take $R\to\infty$. Theorem~\ref{T5.1} controls $Q(V^{(R)})$ and since $\kappa_j^{(R)}$ converge
individually to the $\kappa_j$ associated to $V$, we have
\begin{equation} \lb{7.4}
\sum_{j=1}^\infty \kappa_j^3 \leq  \liminf_{R\to\infty}\, \sum_{j=1}^\infty
\bigl[\kappa_j^{(R)}\bigr]^3
\end{equation}
(a trivial instance of Fatou's Lemma). Thus \eqref{7.3} becomes
\begin{equation} \lb{7.5}
\tfrac18 \int_0^\infty V(x)^2 \, dx \geq \tfrac23 \sum_{j=1}^\infty
\kappa_j^3 + Q(V).
\end{equation}

In particular, $V\in L^2$ implies $Q<\infty$, that is, \eqref{1.21} holds.
As $\kappa_j^3 = [E_j^{(0)}]^{3/2}$, so $\sum [E_j^{(0)}]^{3/2}
<\infty$. By \eqref{1.29}, this implies \eqref{1.20}. \hfill \qed

\subsection*{(i)--(iv) $\bf\Rightarrow V\in L^2$} By Theorem~\ref{T6.1}, $d\rho$ is
the spectral measure of a $V\in L_\loc^2$ so, in particular, \eqref{4.28} holds.
Since $-\kappa_j(t)^3 \leq 0$ and $\log[T(\xi,t)]\leq 0$, this implies that
\begin{equation} \lb{7.6}
\tfrac18 \int_0^t V(x)^2\, dx \leq \tfrac23 \sum_j \kappa_j(0)^3 +
\lim_{y\to\infty} \int_0^\infty P(\xi,y) \log \biggl( \f{1}{T(\xi,0)}\biggr)\,d\xi.
\end{equation}
By the same monotone convergence argument used in the first part of the proof,
\begin{equation} \lb{7.7}
\tfrac18 \int_0^t V(x)^2\, dx \leq \tfrac23 \sum_j [E_j^{(0)}]^{3/2} + Q.
\end{equation}
Taking $t\to\infty$, we see $V\in L^2$ and that \eqref{7.7} holds with $t=\infty$.
\hfill \qed

\smallskip

Our proof shows that the Faddeev--Zhabat sum rule, \eqref{1.27}, holds for any
$V\in L^2 (0,\infty)$. Rewriting $Q$ in terms of the reflection coefficient
(see \eqref{5.13}) and fixed on $(-R,\infty)$ with $R<\infty$, one can obtain
\eqref{1.27} for $V\in L^2 (-\infty,\infty)$ by using the ideas in \cite{SZ}.

\section{Isolating $\Re(w)$} \lb{s8}

The next four sections are devoted to deducing Theorem~\ref{T1.3} from Theorem~\ref{T1.2}.
This amounts to showing that the two lists of conditions are equivalent.  (They are
not equivalent item by item, only collectively.)

The role of this section is to prove that
\begin{equation} \lb{8.2}
\text{Strong Quasi-Szeg\H{o}} \Leftrightarrow \text{Quasi-Szeg\H{o}} + (R<\infty)
\end{equation}
where
\begin{equation} \lb{8.1}
R=\int_0^\infty \log \biggl\{1+ \biggl(\f{\Re w}{k}\biggr)^2\biggr\} k^2\, dk.
\end{equation}
Note that the strong quasi-Szeg\H{o} condition involves both the real and imaginary parts of $w$,
whereas the quasi-Szeg\H{o} condition depends only on $\Im w$ and $R$ only on $\Re w$.  Hence the
title of this section.

We now present an outline of the proof of Theorem~\ref{T1.3}:
under the assumption of the Weyl and Lieb-Thirring conditions, we prove
\begin{SL}
\item Strong Quasi-Szeg\H{o} $\Rightarrow$ Quasi-Szeg\H{o} (this section).
\item Strong Quasi-Szeg\H{o} $+$ Local Solubility $\Rightarrow$ Normalization (see Section~\ref{s11}).
\item Quasi-Szeg\H{o} $+$ Normalization $\Rightarrow R<\infty$ (see Section~\ref{s11}),
so, by \eqref{8.2}, Quasi-Szeg\H{o} $+$ Normalization $\Rightarrow$ Strong Quasi-Szeg\H{o}.
\item Normalization $\Rightarrow$ Local Solubility (see Section~\ref{s9}).
\end{SL}
The first two statements show that the conditions in Theorem~\ref{T1.3} imply those in
Theorem~\ref{T1.2}, the second pair proves the converse.

\begin{lemma}\label{Llog}
For any $f\in\Cmplx$ and any $0<\epsilon\leq 1$,
$$
\log( 1 + |f|^2 ) \leq \begin{cases} \epsilon^{-1} \log( 1 + \epsilon |f|^2 ) & \\
 \log( 1 + \epsilon |f|^2 ) + \log(1 + \epsilon^{-1} ). & \end{cases}
$$
Moreover, if $\epsilon=(1+\delta)^{-2}$ and $\delta \geq 6$, then
$$
\log(1 + \epsilon^{-1} ) \leq 6 \log( \tfrac14\delta + \tfrac12 + \tfrac14\delta^{-1} ).
$$
\end{lemma}

\begin{proof}
The first inequality follows from the concavity of $F:x\mapsto\log(1+x|f|^2)$:
$$
\epsilon \log( 1 + |f|^2 ) = (1-\epsilon) F(0) + \epsilon F(1)
  \leq F(\epsilon) = \log( 1 + \epsilon |f|^2 ).
$$
The second inequality follows from
$$
1 + |f|^2 \leq 1 + \epsilon |f|^2 + \epsilon^{-1} + |f|^2
  = ( 1 + \epsilon |f|^2 )(1 + \epsilon^{-1} )
$$
by taking logarithms.  For the last inequality, notice that
$$
1+(1+\delta)^2 = \delta^2 + 2\delta + 2
    \leq \delta^2 + 4\delta + 6 + 4 \delta^{-1} + \delta^{-2}
    = 16 ( \tfrac14\delta + \tfrac12 + \tfrac14\delta^{-1} )^2
$$
and since $\delta\geq6$, we have $2 \leq \tfrac14\delta + \tfrac12 + \tfrac14\delta^{-1}$.  Therefore,
$$
1+(1+\delta)^2 \leq ( \tfrac14\delta + \tfrac12 + \tfrac14\delta^{-1} )^6,
$$
which gives the result.
\end{proof}

\begin{theorem}\lb{T8.2} Using the notations
\begin{align}
SQS &= \int \log \biggl[ \f{\abs{w(k+i0) + ik}^2}{4k\Im w(k+i0)}\biggr] k^2\, dk,  \\
QS &= \int_0^\infty \log \biggl[ \f{1}{4}\, \f{d\rho}{d\rho_0} + \f12 + \f14\, \f{d\rho_0}{d\rho} \biggr]
        \sqrt{E}\, dE,
\end{align}
and $R$ as in \eqref{8.1}, we have $QS \leq SQS \leq QS + R$ and $R\leq 55 \, SQS$.
In particular, $QS + R <\infty \Leftrightarrow \text{SQS} <\infty$.
\end{theorem}

\begin{proof}
The bulk of the proof rests on the following calculation:
\begin{align}
\f{\abs{w(k + i0) + ik}^2}{4k \Ima w(k+i0)}
&= \f{\Ima w}{4k} + \f12 + \f{k}{4\Ima w} + \f{k}{4\Ima w} \biggl( \f{\Re w}{k}\biggr)^2 \lb{8.3} \\
&= \biggl( \frac{\delta}{4} + \frac12 + \frac{1}{4\delta}\biggr)
    \biggl[ 1 + \frac{1}{(1+\delta)^2}\biggl( \f{\Re w}{k}\biggr)^2 \biggr] \lb{8.4}
\end{align}
where $\delta=\f{\Ima w}{k} = \f{d\rho}{d\rho_0}$.  Taking logarithms and integrating immediately
shows that $QS \leq SQS \leq QS + R$.

To prove $R\leq 55 \, SQS$, we make use of the following notation:
$$
  \delta=\frac{d\rho}{d\rho_0},\quad \epsilon=( 1 + \delta )^{-2}, \quad f(k)=\frac{\Re w}{k},\
  \quad\text{and}\quad A=\{E: \delta > 6 \}.
$$
Notice that from the calculation above,
$$
SQS = \int \log( \tfrac14\delta + \tfrac12 + \tfrac14\delta^{-1}) k^2\,dk
    + \int \log( 1 + \epsilon |f|^2 ) k^2\,dk.
$$
Combining this with Lemma~\ref{Llog} gives
\begin{align*}
R &=\int_0^\infty \log( 1 + |f|^2 ) k^2\,dk \\
    &\leq \int_A \log( 1 + \epsilon |f|^2 ) k^2\,dk + \int_A \log( 1 + \epsilon^{-1} ) k^2\,dk  + {}\\
        &\qquad + 49 \int_{A^c} \log( 1 + \epsilon |f|^2 ) k^2\,dk \\
    &\leq 49 \, SQS + 6 \int_A \log( \tfrac14\delta + \tfrac12 + \tfrac14\delta^{-1}) k^2\,dk \\
    &\leq 55 \, SQS.
\end{align*}
The number $49$ appears because on $A^c$, $\delta\leq6$ which implies $\epsilon^{-1} \leq 49$.
\end{proof}

\section{The Normalization Conditions} \lb{s9}

In this section, we will prove that
$$
 \text{Normalization}\ \Rightarrow \ \text{\eqref{1.17c}}\ \Rightarrow \ \text{Local Solubility}
$$
(cf. step~(iv) in the strategy of Section~\ref{s8}).  This then implies that $d\rho$ is
the spectral measure of a potential $V\in L^2_\loc$ by Theorem~\ref{T6.1}.

\begin{proposition}\lb{P9.1} Let $d\nu$ be any real signed measure on $[0,\infty)$ and
define $M_l\nu$ by \eqref{1.17f}. Then the following are equivalent:
\begin{gather}
M_l\nu \in L^2 (dk) \lb{9.1A} \\
\abs{\nu} ([n,n+1])\in \ell^2 \lb{9.1B} \\
\int \log \biggl[ 1+\biggl( \f{M_l\nu}{k}\biggr)^2\biggr] k^2\, dk <\infty. \lb{9.1C}
\end{gather}
\end{proposition}

\begin{proof}  It is not difficult to see that  \eqref{9.1A} $\Rightarrow$ \eqref{9.1B}:
\begin{align}\label{9.4}
\sum   \bigl[ \abs{\nu} ([n,n+1]) \bigr]^2 &\leq \int_0^\infty   \bigl[ \abs{\nu} ([k-1,k+1]) \bigr]^2 \, dk
\leq 4 \int_0^\infty   \bigl[ M_l\nu(k) \bigr]^2 \, dk.
\end{align}
To prove the converse, let us write $\nu_n =\abs{\nu} ([n,n+1])$.  Then, for any $k\in[n,n+1]$,
\begin{equation}\label{E:DMO}
M_l\nu(k) = \sup_{L\geq 1} \frac{|\nu|([k-L,k+L])}{2L}
    \leq \sup_{m\geq 0} \frac{3}{2m+1} \sum_{j=n-m}^{n+m} \nu_j.
\end{equation}
Indeed, one may take $m$ to be the integer in $[L,L+1)$.  As the
discrete maximal operator in \eqref{E:DMO} is $\ell^2$ bounded, we may deduce
\begin{align}
\int_0^\infty \bigl[ M_l\nu(k) \bigr]^2 \, dk \leq \sum_{n} \Bigl[ \sup_{k\in[n,n+1]} M_l\nu(k) \Bigr]^2
&\leq C \sum_{n} \nu_n^2.
\end{align}
This proves \eqref{9.1B} $\Rightarrow$ \eqref{9.1A}.

As $\log(1+x^2)\leq x^2$,
$$
\int \log[ 1 + k^{-2} M_l\nu(k)^2 ] \,k^2\, dk  \leq \int [ M_l\nu(k) ]^2 \, dk,
$$
which proves \eqref{9.1A} $\Rightarrow$ \eqref{9.1C}.

We will finish the proof by showing that \eqref{9.1C} implies \eqref{9.1B}.  For each $k\in[n,n+1]$,
it follows directly from the definition that $\tfrac12 \nu_n \leq M_l\nu(k)$.  Thus
\begin{equation}\label{E:9.7}
\sum_n  n^2 \log\biggl[ 1 + \frac{\nu_n^2}{4(n+1)^2} \biggr]
\leq  \int \log \biggl[ 1+\biggl( \f{M_l\nu}{k}\biggr)^2\biggr] k^2\, dk,
\end{equation}
which shows that $\nu_n \geq (n+1)$ only finitely many times.  For the remaining values of $n$, one need
only apply the estimate $\log(1+x)\geq \tfrac12 x$ for $x\in[0,1]$, which follows by comparing derivatives,
to see that $\nu_n\in\ell^2$.
\end{proof}

\begin{theorem}\lb{T9.2} If
\[
\int \log\biggl[ 1+\biggl(\f{M_s\nu}{k}\biggr)^2\, \biggr] k^2\, dk <\infty,
\]
then the equivalent conditions of Proposition~\ref{P9.1} hold.
\end{theorem}

\begin{proof} The result follows by the reasoning used to prove \eqref{9.1C} $\Rightarrow$ \eqref{9.1B}:
For all $k\in[n,n+1]$,
$$
\abs{\nu} ([n,n+1])  \leq  \abs{\nu} ([k-1,k+1]) \leq 2 M_s\nu(k).
$$
Thus \eqref{E:9.7} holds with $M_s\nu$ in place of $M_l\nu$ and the argument given above may be continued from there.
\end{proof}

\begin{theorem}\lb{T9.3} If \eqref{9.1B} holds, then so does Local Solubility, that is,
\eqref{1.19x}. In particular, by Theorem~\ref{T9.2},
\[
\text{Normalization} \Rightarrow \text{Local Solubility}.
\]
\end{theorem}

\begin{remark} This is step (iv) of the strategy in Section~\ref{s8}.
\end{remark}

\begin{proof} By the definition \eqref{1.18},
\begin{equation}\lb{9.21x}
  F(q) =2\pi^{-1/2} \int_{p\geq 1} e^{-(p-q)^2} \, d\nu(p).
\end{equation}
As $(n+x-q)^2 \geq (n-q)^2 -  2\abs{x}\abs{n-q}$ for $|x|<1$,
\begin{equation} \lb{9.21}
\abs{F(q)}\leq 2\pi^{-1/2} \sum_{n=1}^\infty e^{-(n-q)^2} e^{2\abs{n-q}}
\abs{\nu} ([n,n+1]).
\end{equation}
Thus by Young's inequality for sums, \eqref{9.1B} $\Rightarrow F\in L^2$.
\end{proof}

We conclude this section with a result we will need in Section~\ref{s11}.  In the proof,
we will use the following simple inequality: for $\delta\in[0,1]$,
\begin{equation}\label{E:SI}
\log [\tfrac14\, \delta+\tfrac12 + \tfrac14\, \delta^{-1}] \geq \tfrac14(\delta-1)^2.
\end{equation}
As equality holds when $\delta=1$, the result follows by differentiating:
$$
\frac{\delta-1}{\delta(\delta+1)} \leq \frac12(\delta-1).
$$

\begin{theorem}\lb{T9.4} If {\rm{(}}Strong{\rm{)}} Quasi-Szeg\H{o} and Local Solubility hold,
then $\nu$ obeys \eqref{9.1B}. In particular, by Theorem~\ref{T1.3}, this follows for $V\in L^2$.
\end{theorem}

\begin{proof} Let us recall that the Quasi-Szeg\H{o} condition says
\begin{equation}\label{E:QS9}
\int_0^\infty \log \biggl[ \f{1}{4}\, \f{d\rho}{d\rho_0} + \f12 + \f14\, \f{d\rho_0}{d\rho} \biggr] k^2\,dk < \infty.
\end{equation}
(By Theorem~\ref{T8.2}, this is also implied by the strong quasi-Szeg\H{o} condition.)

Let us decompose $d\nu=d\nu_+ - d\nu_-$ where $d\nu_\pm$ are both positive measures.
The definition of $d\nu$, \eqref{1.17a}, shows that for $k>1$,
$$
\frac{d\rho}{d\rho_0} = 1 + \frac{1}{k}\frac{d\nu}{dk}.
$$
Moreover, $d\nu_-$ is absolutely continuous; in fact, \eqref{1.17b} shows $\frac{d\nu_-}{dk}\leq k$.

Let us restrict the integral \eqref{E:QS9} to the essential support of $d\nu_-$,
that is, where $\f{d\rho}{d\rho_0}\leq 1$.  Using \eqref{E:SI}, we deduce that
\begin{equation}
\int_0^\infty   \biggl| \f{d\nu_-}{dk} \biggr|^2 \,dk < \infty
\end{equation}
and hence that $\abs{\nu_-} ([n,n+1]) \in \ell^2$.  To complete the proof, we need to deduce the same result
for $\nu_+$.

The local solubility condition says $F\in L^2$ where $F$ is defined as in \eqref{9.21x}.
The first sentence of Theorem~\ref{T9.3} says
\[
F_-(q) = \int_{p\geq 1} e^{-(p-q)^2}\, d\nu_- (p) \ \in L^2 (dq)
\]
and so $F\in L^2$ implies
\[
F_+(q)=\int_{p\geq 1} e^{-(p-q)^2} \, d\nu_+ (p) \ \in L^2 (dq).
\]
For $q\in[n,n+1]$, we have $F_+(q) \geq e^{-1} \nu_+ ([n,n+1])$ and thus may conclude $\nu_+ ([n,n+1])\in\ell^2$.
\end{proof}

\section{Harmonic Analysis Preliminaries} \lb{s10}

For harmonic functions in the half-plane, it is well known that the conjugate function
belongs to $L^p$ ($0<p<\infty$) if and only if the same is true for the nontangential maximal function.
The first direction appears already in the paper of Hardy and Littlewood that introduced the
maximal function \cite[Theorem 27]{HL}.  The other direction, which is much harder, is due
to Burkholder, Gundy, and Silverstein \cite{BGS}.
The purpose of this section is to present an analogous theorem with a peculiar replacement for $L^p$.
Theorem~2 of \cite{BGS} covers this situation perfectly if one is willing to consider the
maximal Hilbert transform; we are not.  However, this does resolve one direction; for the other, we will
use subharmonic functions in the manner of \cite{HL}.

We will use the following notation: $f\lesssim g$ means $f\leq C g$ for some absolute constant $C$,
whereas $f\approx g$ means that $f\lesssim g$ and $g \lesssim f$.

\begin{proposition}
Let $d\sigma$ be a compactly supported positive measure on $\Reals$,
\begin{equation}\label{BGS1}
    \int \log \bigl[1 + |H\sigma|^2 \bigr] \,dx \lesssim
    \int \log \bigl[1 + |M\sigma|^2 \bigr] \,dx.
\end{equation}
\end{proposition}

\begin{proof}
This is a special case of \cite[Theorem 2]{BGS}.  It is also amenable to the good-$\lambda$
approach discussed in textbooks: \cite[\S V.4]{Stein} or \cite[\S XIII]{Tor}.
\end{proof}

As noted earlier, Burkholder, Gundy, and Silverstein do not provide the converse inequality;
indeed as they note, in the generality they treat, the result is false without switching to the
maximal Hilbert transform.  Nevertheless,
the function $x\mapsto\log[1+x^2]$ grows sufficiently quickly that the result is true.  We
divide the proof into two propositions.

\begin{proposition}\label{P7.3}
There is a $\lambda_0$ so that for any finite positive measure $d\sigma$ on $\Reals$,
\begin{align}\label{Ineq2}
\int_{\{M \sigma > \lambda_0 \}} \log\bigl[ 1 + |M\sigma|^2 \bigr]\,dx \lesssim \int \log \bigl[1 +
\left(H \sigma\right)^2 + \left(\tfrac{d \sigma}{dx}\right)^2 \bigr] \,dx.
\end{align}
In particular, $|\{M \sigma > 2\lambda_0 \}|\lesssim \mbox{\rm RHS\eqref{Ineq2}}$.
\end{proposition}

\begin{proof}
Let $u(z)+iv(z)=\int d\sigma(x)/(x-z)$ denote the Cauchy integral of $d\sigma$, then
$$
  F(z) = \log[1+u(z)+iv(z)]
$$
is analytic---$u\geq 0$ because it is the Poisson integral of a positive measure. In particular,
$|F|^{1/2}$ is subharmonic.  Now as $|F(z)| \geq \log|1+u(z)|$,
\begin{align}
\log\bigl[ 1+[M\sigma](x)\bigr] &\lesssim \sup_{y>0} \log\bigl[ 1+ u(x+iy) \bigr] \\
    &\leq \sup_{y>0} |F(x+iy)| \\
    &\lesssim \bigl\{ [ M |F|^{1/2} ](x) \bigr\}^2.
\end{align}
Elementary calculations show $\abs{\Re F}\leq \log(1+u+|v|)$ and $\abs{\Im F}\leq\frac\pi2$; therefore,
$$
\log\bigl[ 1+ M\sigma \bigr] \lesssim \left\{ 1 + M \sqrt{\log[1+u+|v|]} \right\}^2
\lesssim 1 + \left\{ M \sqrt{\log[1+u+|v|]} \right\}^2.
$$
>From this, one may deduce that for $\lambda_1$ sufficiently large,
\begin{align}\label{E:star}
\sqrt{\log\bigl[ 1+ M\sigma \bigr]} \lesssim M \sqrt{\log[1+u+|v|]}
\end{align}
on the set where $\log[ 1+M\sigma]\geq \lambda_1$.

Interpolating between the $L^\infty$ and $L^2$ bounds on $M$ shows that
$$
\int_{Mf>\lambda} |M f|^2 \,dx  \lesssim \int_{|f|>\lambda/2} |f|^2 \, dx.
$$
Combining this with \eqref{E:star}, we see that for $\lambda_0\geq e^{\lambda_1}$ and $\epsilon$
sufficiently small,
$$
    \int_{\{M \sigma > \lambda_0 \}} \log\bigl[ 1+ M\sigma \bigr] \,dx
\lesssim \int_{\{ |u+iv|>\epsilon \}}  \log[1+u+|v|] \, dx.
$$
To obtain \eqref{Ineq2}, we need merely note that $\log(1+x)\approx\log(1+x^2)$ on any interval
$[a,\infty)$ with $a>0$.
\end{proof}

\begin{proposition}\label{HA.4}
For any finite positive measure $d\sigma$ on $\Reals$,
\begin{equation}
\int \log \bigl[1 + \left(M\sigma\right)^2 \bigr] \,dx \lesssim \int \log \bigl[1 + \left(H
\sigma\right)^2 + \left(\tfrac{d \sigma}{dx}\right)^2 \bigr] \,dx.
\end{equation}
\end{proposition}

\begin{proof}
By Proposition~\ref{P7.3}, it suffices to prove
\begin{equation}\label{want}
\int_{\{M \sigma \leq \lambda_0 \}} \abs{M\sigma}^2 \,dx \lesssim \int \log \bigl[1 + (H \sigma)^2 +
(\tfrac{d \sigma}{dx})^2 \bigr] \,dx.
\end{equation}

Let $\Omega=\{M \sigma > 4\lambda_0 \}$, $d\sigma_1=\chi_\Omega d\sigma$, and
$d\sigma_2=\chi_{\Omega^c}d\sigma$. We will prove \eqref{want} by writing $M\sigma\leq M\sigma_1 +
M\sigma_2$.

It is a well-known property of the maximal function that
$$
\sigma(\{M \sigma > 4\lambda_0 \})\lesssim \lambda_0|\{M \sigma > 2\lambda_0 \}|.
$$
Combining this with Proposition~\ref{P7.3}, shows that $\|\sigma_1\|=\sigma(\Omega)\lesssim
\text{RHS\eqref{want}}$. Consequently, by the weak-type $L^1$ bound on the maximal operator,
$$
\int_{\{M \sigma \leq \lambda_0 \}} \abs{M\sigma_1}^2 \,dx \lesssim \int_0^{\lambda_0}
 \frac{\|\sigma_1\|}{\lambda} 2\lambda\,d\lambda
\lesssim \text{RHS\eqref{want}}.
$$

Now we turn to bounding $M\sigma_2$.  On $\Omega^c$, we know that $d\sigma$ must be absolutely
continuous and its Radon-Nikodym derivative is bounded by $4\lambda_0$.  Therefore, $L^2$
boundedness of the maximal operator implies
$$
 \int \abs{M\sigma_2}^2 \,dx
    \lesssim \int_{\{M\sigma\leq 4\lambda_0\}} |\tfrac{d \sigma}{dx}|^2 \,dx
    \lesssim \int \log \bigl[1 + \left(\tfrac{d \sigma}{dx}\right)^2 \bigr] \,dx,
$$
which completes the proof.
\end{proof}

Putting the previous propositions together, we obtain the following

\begin{theorem}\label{Tloc}
If $\sigma$ is a positive measure of compact support, then
\begin{equation}\label{BGS2}
    \int \log \Bigl[1 + |H\sigma|^2 + (\tfrac{d \sigma}{dx})^2 \Bigr] \,dx \approx
    \int \log \Bigl[1 + |M\sigma|^2 \Bigr] \,dx.
\end{equation}
\end{theorem}

\section{Taming $\Real m$} \lb{s11}

The purpose of this section is to prove Corollary~\ref{C11.3} below and so complete the
proof of Theorem~\ref{T1.2} as laid out in Section~\ref{s8}.

Let $H_s$ denote the short-range Hilbert transform: $H_s \sigma=K*\sigma$ where
$$
K(x) = \begin{cases} 0, & |x|>1 \\ \tfrac1\pi[x^{-1} - x], & |x|<1 \end{cases}
$$
and let $H_l=H-H_s$ denote the long-range Hilbert transform: $H_l\sigma=K*\sigma$ with
$$
K(x) = \begin{cases} \tfrac1\pi x^{-1}, & |x|>1 \\ \tfrac1\pi x, & |x|<1. \end{cases}
$$
Note that both $H_s$ and $H_l$ are Calder\'on--Zygmund operators and so bounded on $L^p(\Reals)$ for
$1<p<\infty$.  As in the Introduction, we define short- and long-range maximal operators:
$$
[M_s \sigma](x) = \sup_{L\leq1} \frac{ |\sigma|\bigl([x-L,x+L]\bigr) }{2L},
$$
and for $M_l$, the supremum is taken over $L\geq 1$.  Naturally, both truncated maximal operators
are $L^p$-bounded for $1<p\leq\infty$.

We will use the notation
$$
    \|\mu\|_{\ell^2(M)}^2 =  \sum_n \Bigl[ |\mu|\bigl( [n,n+1] \bigr) \Bigr]^2
$$
as introduced in \eqref{1.17c}. Obviously, $\|\mu\|_{\ell^2(M)}^2 \leq \|\mu\|^2$.

\begin{lemma}\label{dumb}
Let $F(k) = (1+k^2)^{-1}$.  For each complex measure $\mu\in\ell^2(M)$,
\begin{align*}
\int \bigl[\Phi*|d\mu|\bigr]^2 \,dk \lesssim \|\mu\|_{\ell^2(M)}^2, \ %
\int |M_l\mu|^2 \,dk \lesssim \|\mu\|_{\ell^2(M)}^2, \ %
\int |H_l\mu|^2 \,dk \lesssim \|\mu\|_{\ell^2(M)}^2.
\end{align*}
\end{lemma}

\begin{proof}
All three inequalities follow by replacing $|d\mu|$ by its average on each of the intervals $[n,n+1]$.
This operation changes $\Phi*|d\mu|$
and $M_l\mu$ by no more than a factor of two.  For $H_l$, it introduces an error which can be
bounded by $\Phi*|d\mu|$.  We then use the $L^2$ boundedness of the appropriate operator.
\end{proof}

\begin{theorem}\label{T:IT}
If $\mu$ is a positive measure on $\Reals$ with $\|\mu\|_{\ell^2(M)}^2<\infty$, then
\begin{equation}\label{BGSa}
    \int \log \Bigl[1 + |H\mu|^2k^{-2} \Bigr] \,(1+k^2)\,dk < \infty
\end{equation}
if and only if
\begin{equation}\label{BGSb}
    \int \log \Bigl[1 + |M_s\mu|^2 k^{-2} \Bigr] \,(1+k^2)\,dk < \infty.
\end{equation}
\end{theorem}

\begin{proof}
As neither integral can diverge on any compact set, we can restrict our attention to $k>1$.

We begin by proving that \eqref{BGSb} implies \eqref{BGSa}.  Given a compactly supported positive measure $d\sigma$,
Theorem~\ref{Tloc} and Lemma~\ref{dumb} show that
\begin{align*}
 \int_n^{n+1} \log  \Bigl[1 + |H_s\sigma|^2 \Bigr] \,dk
&\lesssim \int_n^{n+1} \log \Bigl[1 + |H\sigma|^2 \Bigr] + |H_l\sigma|^2 \,dk  \\
&\lesssim \|\sigma\|^2_{\ell^2(M)} + \int_n^{n+1} \log \Bigl[1 + |H\sigma|^2 \Bigr] \,dk \\
&\lesssim \|\sigma\|^2_{\ell^2(M)} + \int \log \Bigl[1 + |M\sigma|^2 \Bigr] \,dk \\
&\lesssim \|\sigma\|^2_{\ell^2(M)} + \int \log \Bigl[1 + |M_s\sigma|^2 \Bigr] \,dk + \int |M_l\sigma|^2 \,dk \\
&\lesssim 2\|\sigma\|^2_{\ell^2(M)} + \int \log \Bigl[1 + |M_s\sigma|^2 \Bigr] \,dk.
\end{align*}
Let us choose $\sigma=(1+n^2)^{-1/2}d\mu_n$ where $d\mu_n$ is the restriction of $d\mu$ to the interval
$[n-1,n+2]$.  Combining the above with Lemma~\ref{dumb} gives
\begin{align*}
\sum (1+n^2) & \int_n^{n+1} \log \Bigl[1 + \tfrac1{n^2+1}|H\mu|^2 \Bigr] \,dk  \\
&\lesssim \|\mu\|_{\ell^2(M)}^2 + \sum (1+n^2) \int_n^{n+1} \log \Bigl[1 + \tfrac1{n^2+1}|H_s\mu|^2 \Bigr] \,dk \\
&\lesssim \|\mu\|_{\ell^2(M)}^2 + \sum (1+n^2) \int_{n-3}^{n+4} \log \Bigl[1 + \tfrac1{n^2+1}|M_s\mu|^2 \Bigr] \,dk \\
&\lesssim \|\mu\|_{\ell^2(M)}^2 + \int \log \Bigl[1 + k^{-2} |M_s\mu|^2 \Bigr] (1+k^2)\,dk.
\end{align*}

The proof that \eqref{BGSa} implies \eqref{BGSb} is a little more involved because the Hilbert transform is
positivity preserving.

Let $\phi$ be a smooth bump which is supported on $[-2,3]$ and is equal to $1$ on $[-1,2]$.  We will write
$\phi_n(x)$ for $\phi(x-n)$.  Elementary calculations show that
\begin{equation}\label{ssss}
\Bigl| [H_s(\phi d\sigma)](x) - \phi(x)[H_s \sigma](x) \Bigr| \lesssim \|\sigma \| \Phi(x)
\end{equation}
where $\Phi(x)=(1+x^2)^{-1}$. Using Theorem~\ref{Tloc}, Lemma~\ref{dumb}, and then \eqref{ssss},
\begin{align*}
\int_n^{n+1} \log \Bigl[1 + |M_s\sigma|^2 \Bigr] \,dk
&\leq \int_n^{n+1} \log \Bigl[1 +  |M(\phi_nd\sigma)|^2 \Bigr] \,dk \\
&\lesssim \int \log \Bigl[1 + |H(\phi_nd\sigma)|^2 + |\phi_n\tfrac{d\sigma}{dk}|^2 \Bigr] \,dk \\
&\lesssim \int \log \Bigl[1 + |H_s(\phi_nd\sigma)|^2 \Bigr] \,dk +  \|\sigma\|^2\\
&\lesssim \int \log \Bigl[1 + \phi_n^2 |H_s\sigma|^2 \Bigr] \,dk + \|\sigma\|^2
\end{align*}
By choosing $\sigma=(1+n^2)^{-1/2}d\mu_n$ where $d\mu_n$ is the restriction of $d\mu$ to the interval
$[n-4,n+5]$, the proof may be completed in much the same manner as was used to prove the opposite
implication.
\end{proof}

It is now easy to complete the outline from Section~\ref{s8}.

\begin{corollary}\label{C11.3} In the nomenclature of Theorems~\ref{T1.2}, \ref{T1.3}, and~\ref{T8.2},
\begin{align}
\label{E:C1}  \text{Normalization} &\Rightarrow R<\infty \\
\label{E:C2} \text{Strong Quasi-Szeg\H{o}} \ + \  \text{Local Solubility} &\Rightarrow \text{Normalization}.
\end{align}
\end{corollary}

\begin{proof}
We begin with \eqref{E:C1}.  As the $m$-function associated
to the free operator is purely imaginary on the spectrum, we have that for all $k>0$,
\begin{align} \label{E:R1}
\Re w(k) &= \int_{(-\infty,1]} \frac{d\rho(E)-d\rho_0(E)}{E-k^2}  + \frac2\pi \int \frac{\xi d\nu(\xi)}{\xi^2-k^2} \\
&= f(k) + \frac1\pi \int \frac{d\nu(\xi)}{\xi+k} + \frac1\pi \int \frac{d\nu(\xi)}{\xi-k} \\
&= f(k) + [H\mu](k)
\end{align}
where $f(k)$ is defined to be the first term on the RHS of \eqref{E:R1} and $d\mu$ is defined by
$$
\int g(k) \,d\mu(k) = \int g(k) \,d\nu(k) + \int g(-k) \,d\nu(k).
$$

By Theorem~\ref{T9.2}, Normalization implies $\nu\in\ell^2(M)$ and hence $\mu\in\ell^2(M)$; thus
we may apply Theorem~\ref{T:IT} to see that \eqref{BGSa} holds.  As
$$
\log[1+(x+y)^2] \lesssim x^2 + \log[1+y^2]
$$
and $|f(k)|\lesssim(k-1)^{-1}$ for $k>1$, we see that this is sufficient to deduce $R<\infty$.

We now turn to \eqref{E:C2}.  By Theorem~\ref{T8.2}, we know that $R<\infty$ and so by the calculation above,
\eqref{BGSa} holds.  From the proof of Theorem~\ref{T9.4} we are guaranteed that $d\mu$, defined as
above, belongs to $\ell^2(M)$.  Thus we may apply Theorem~\ref{T:IT} to deduce that the normalization condition
holds.
\end{proof}


\end{document}